\documentclass[10pt,twoside]{article}

\usepackage{amsmath}
\usepackage{amsfonts}
\usepackage{amssymb}
\usepackage{amscd}
\usepackage{amsthm}
\usepackage{amsbsy}
\usepackage{graphicx}
\usepackage{bm}

\def\@oddhead{\hfill \shorttitle \hfill \thepage}
\def\@evenhead{\thepage \hfill \shortauthor \hfill}
\def\@oddfoot{}
\def\@evenfoot{}

\newtheorem{theo}{Theorem}[section]

\newtheorem{teo}[theo]{Theorem}
\newtheorem{lem}[theo]{Lemma}

\theoremstyle{definition}
\newtheorem{defi}[theo]{Definition}

\newtheorem{rema}[theo]{Remark}

\newtheorem{remark}[theo]{Remark}

\newcommand{\ADS}{\mathbb{ADS}}
\newcommand{\DS}{\mathbb{DS}}

\newcommand{\wt}{\widetilde}

\newcommand{\cT}{\mathcal T}

\newcommand{\DD}{\mathbb{D}}

\newcommand{\HH}{\mathbb{H}}

\newcommand{\RR}{\mathbb{R}}
\renewcommand{\SS}{\mathbb{S}}

\newcommand{\ZZ}{\mathbb{Z}}
\newcommand{\op}{\operatorname}

\DeclareMathOperator{\SO}{SO}

\DeclareMathOperator{\AdS}{AdS}
\DeclareMathOperator{\dS}{dS}

\DeclareMathOperator{\Ein}{Ein}

\newcommand{\ie}{\textit{i.e. }}

\newcommand{\sCa}[2]{\langle #1 \;|\; #2\rangle}         

\newcommand{\mcal}[1]{\ensuremath{\mathcal{#1}}}
\newcommand{\orth}{\bot}

\newcommand{\cD}{\mcal{D}}

\newcommand{\cK}{\mcal{K}}

\newcommand{\cC}{\mcal{C}}
\newcommand{\cH}{\mcal{H}}

\date{}
\title{\ \\[0.4cm] \ \\ \bf  Lorentzian Kleinian Groups }
\author{Thierry Barbot\footnote{LMA,
Avignon University,  Campus Jean-Henri Fabre
301 rue Baruch de Spinoza
BP 21239
84916 Avignon Cedex 9,
France. E-mail: Thierry.Barbot@univ-avignon.fr.
}}
\begin{document}

\maketitle


\thispagestyle{empty}

\begin{abstract}
\vskip 3mm\footnotesize{

\vskip 4.5mm
\noindent Classical Kleinian groups are discrete subgroups of isometries of $\HH^n$. The well-known theory of Kleinian groups starts with the definition of their associated limit set in the boundary of $\HH^n$, and includes the geometric properties of the quotient hyperbolic space.

This approach, naively applied, fails in the Lorentzian analogue anti-de Sitter space: discrete subgroups
do not act properly discontinuously, and in many cases the set of accumulation points of orbits at the conformal boundary at infinity depends on the orbit.

In this survey, we point out a way to extend this classical theory by introducing causality notions: the theory of limit sets and regularity domains extend naturally to achronal subgroups. This is closely related to the notions of globally hyperbolic spacetimes, and we present what is known about the classification of globally hyperbolic spacetimes of constant curvature. We also review the close connection revealed by G. Mess (\cite{mess1}) between globally hyperbolic spacetimes of dimension $2+1$ and Teichm\"uller space. This link can be understood via the space of timelike geodesics of anti-de Sitter space, and this space has also an interesting role, presented here, in the recent works about proper group actions on spacetimes of constant curvature.

\vspace*{2mm} \noindent\textbf{\bf 2000 Mathematics Subject Classification:} 30F40, 30F60, 53B30, 53C25, 53C35, 53C50, 53C55, 53D30, 57M50, 83C05.

\vspace*{2mm} \noindent\textbf{\bf Keywords and Phrases: }Discrete groups, Lorentzian spacetimes, causality notions, global hyperbolicity, Anosov representations, Teichm\"uller space.}

\end{abstract}

\tableofcontents
\section{Introduction}
\label{sec:intro}
Kleinian groups\index{Kleinian group}, \ie discrete groups of isometries of the hyperbolic space\index{hyperbolic space} $\HH^n$, are central objects of study in geometry, most often in low dimension, and particularly in dimension $3$. They have several ramifications in other fields of mathematics: number theory, topology, and, of course, group theory.
For nice surveys on this very rich topic, including reports on recent results in this field, let us mention the references \cite{martin, ohshika} in this series of handbooks dedicated to group actions.

The isometry group of $\HH^n$ is $\SO(1,n)$, and here we consider the analogous case of discrete subgroups of
$\SO(2,n)$, that we call in this survey \emph{Lorentzian Kleinian groups\index{lorentzian Kleinian group}.} The group $\SO(2,n)$ is naturally the group of isometries of the anti-de Sitter space, denoted by $\AdS^{1,n}$. Very little is known on Lorentzian Kleinian groups, if compared with the venerable hyperbolic case. The purpose of this survey is to present
part of this relatively young matter.

Classical Kleinian groups have several nice basic properties:

-- the action on $\HH^n$ is properly discontinuous;

-- for any point $x$ in $\HH^n$, the orbit of $x$ under the action of the group accumulates at the conformal boundary $\partial\HH^n \approx \SS^{n-1}$ on a closed invariant set, called the limit set, which does not depend on $x$.

These properties completely fail in the AdS case: stabilizers of points are not anymore compact. Therefore, the action is proper only in particular cases, and the asymptotic behavior of orbits is not uniform: for example a Lorentzian Kleinian group may have infinite orbits and other orbits reduced to a point.

Hence, at first glance, the Lorentzian case seems radically different, and some geometers may consider it as a quite completely different field. As a matter of fact, the first works on Lorentzian geometry in the case of constant curvature were mostly devoted to the study of \emph{compact} Lorentzian manifolds, which has very few similarity with the Riemannian case. Even the geodesic completeness of compact Lorentzian manifolds of constant curvature, even if true, is far from being trivial (\cite{carriere, klingler}) --- and as a matter of fact, compact Lorentzian manifolds are \emph{not} geodesically complete in general.
Wolf's book has long been --- and still is! --- one of the main references in the field, and was mainly concerned with the compact case. One indication on the Riemannian oriented spirit of Wolf's book \cite{wolf} is that it is also one of the main references for the classification of \emph{Riemannian} crystallographic groups.

Nevertheless, there are several important common features between the Riemannian and the Lorentzian case.

-- anti-de Sitter space admits a natural conformal boundary: the Einstein universe $\Ein^{1,n-1}$;

-- in the same way as the Euclidean space and the hyperbolic space embed conformally in the sphere, the anti-de Sitter space embeds in the Einstein space, and moreover, the Minkowski space and the de Sitter space, Lorentzian analogues of the Euclidean space and the sphere, embed in $\Ein^{1,n}$ as well.

The spirit underlying this survey is that there is a fundamental framework in which Riemannian or Lorentzian Kleinian groups have the same nature, but involving a notion that is unapparent because trivial in the
Riemannian case: the notion of \emph{causality.}

It is good in this celebration's year of General Relativity\footnote{This survey has been written in 2015.} to recall that Loren\-tzian geometry is the geometry of space \emph{and time}; in which the classical Riemannian geometry is enclosed as the  \emph{static} case, \ie the case in which the space does not change with time.

This point of view, developped in Section \ref{sec:discrete}, is that one should distinguish certain subgroups, the \emph{achronal subgroups,} that have a reasonable behavior relatively to causality.
Riemannian Kleinian groups are automatically achronal, whereas Lorentzian ones may not be --- as a matter of fact, Lorentzian Kleinian groups acting cocompactly are \emph{never} achronal. Achronal Kleinian groups then appear as completely similar to their Riemannian counterparts: they admit a limit set in the conformal boundary, and they do act properly, not on the entire space itself, but on a certain domain: the domain of points invisible from the limit set.

Therefore, we start this paper by a review on the causality notions (Section \ref{sec:introcausal}), followed by a presentation of the Lorentzian spacetimes of constant curvature (Section \ref{sec:models}).
What appear as the true analogues of compact Riemannian spaces of constant curvature are not the compact Lorentzian manifolds, but the \emph{maximal globally hyperbolic spatially compact spacetimes} (abbreviation MGHC), a notion arising from the analytic treatment of General Relativity with the tools of Partial Differential Equations. In Section \ref{sec:classGH}, we present the classification of MGHC spacetimes of
constant curvature. This classification, initiated by the fundamental work of G. Mess (\cite{mess1, mess2}), is quite recent, and still incomplete. It provides an interesting framework in which many results or interrogations concerning hyperbolic Kleinian groups may very well find in a near future a natural continuation and extension.

After Section \ref{sec:discrete} mentioned above devoted to achronal Kleinian groups (but where we also mention links with the theory of Anosov representations), Section \ref{sec:trois} points out very interesting connections, initiated by the pioneering work of Mess, between MGHC spacetimes of constant curvature and the Teichm\"uller space\index{Teichm\"uller space} Teich$(S)$ of a closed surface $S$: it generalizes the link between hyperbolic $3$-manifolds and Teich$(S)$, appearing for example in the case of quasi-fuchsian manifolds. But it also provides new ones, to which we propose a brief introduction. One explanation behind this phenomenon is the fact that $\SO(2,n)$ has rank $2$, whereas $\SO(1,n)$ has ``only'' rank $1$. In other words, unlike $\SO(1,n)$, $\SO(2,n)$ has not only one, but two maximal parabolic subgroups, meaning that there is not only one geometry, like the hyperbolic geometry, associated with $\SO(2,n)$, but \emph{two} geometries: anti-de Sitter geometry, and also the geometry of the symmetric space $\cT_{2n}$ associated with $\SO(2,n)$, which happens to be \emph{the space of AdS timelike geodesics.} For $n=2$, $\cT_4$ is nothing but the product $\HH^2 \times \HH^2$. The space $\cT_{2n}$ will appear in
this survey in two situations:

--- there is an interplay between the AdS MGHC spacetime $M_\Gamma$ associated with a Lorentzian Kleinian group $\Gamma$ and the quotient $\Gamma\backslash\cT_{2n}$, emerging through Cauchy hypersurfaces, \ie isometric embeddings of Riemannian hypersurfaces in $M_\Gamma$. More precisely, $\Gamma\backslash\cT_{2n}$ has a natural K\"ahler structure, and there is an almost correspondence between Lagrangian submanifolds in $\Gamma\backslash\cT_{2n}$ and Cauchy hypersurfaces in $M_\Gamma$. Moreover, in the case $n=2$,
the symplectic form on $\Gamma\backslash\cT_{2n} \approx \Gamma\backslash(\HH^2 \times \HH^2)$ is the difference $p_1^*\omega_0 - p_2^*\omega_0$ where $\omega_0$ is the volume form on $\HH^2$, and $p_1$, $p_2$ the projections on the left and right factors, respectively. Therefore,
typical Lagrangian submanifolds of $\Gamma\backslash\cT_{2n} \approx \Gamma\backslash(\HH^2 \times \HH^2)$ are graphs
of volume preserving maps between hyperbolic surfaces.
This connection between AdS geometry and special volume preserving maps between hyperbolic surfaces will be developed further, but with a quite different point of view,  in \cite{fillgraham}, to appear in one of the \emph{Handbooks of group actions.}

--- $\cT_{2n}$ is also related to questions of proper actions: at least in dimension $2+1$, a Lorentzian Kleinian group $\Gamma$ acts properly discontinuously on $\AdS^{1,2}$ if and only if there
is an embedded surface $S$ in $\Gamma\backslash(\HH^2 \times \HH^2)$ such that the restriction
to $S$ of the pseudo-Riemannian metric $g_{hyp} - g_{hyp}$ is positive definite - where $g_{hyp}$ denotes the hyperbolic metric, and $g_{hyp} - g_{hyp}$ a simplified notation for $p_1^*g_{hyp} - p_2^*g_{hyp}$. This criterion is related to the existence of $\Gamma$-invariant foliations by timelike geodesics, and we conjecture an extension of this fact in higher dimensions (see Remark \ref{rk:zeghibfoliate}).

We conclude the survey with a quick overview on Lorentzian Kleinian groups acting properly, even cocompactly. This topic has received recent brilliant contributions by J. Danciger, F. Gu\'eritaud and F. Kassel that we mention very briefly, as important illustrations of the role of the space of timelike geodesics. There are several more or less recent surveys on this topic in which the reader may find more substantial information  (\cite{surveyBZ, bigsurvey, schbourbaki}).

The topic is currently growing quite quickly, and this survey has no pretention to be complete. We would like to attract the attention to a very recent work, transposing to Anti-de Sitter space the notion of Patterson-Sullivan measures, and establishing in this context an inequality between the critical exponent and the Hausdorff dimension of the (acausal) limit set (\cite{glorieux}).
\medskip

\noindent\textbf{Remarks on conventions and notation.}
A {\it Lorentzian manifold}\index{Lorentzian manifold} is a manifold equipped with a pseudo-Riemannian metric of signature $(1,n)$ for
some $n\geq 1$. The manifold is then of dimension $n+1$: $n$-dimensional in space, and $1$-dimensional in time.
In our convention a Lorentzian metric has signature $(-, +, ... ,+)$; the value of the metric on a tangent vector $v$ is called the \textit{norm} of
$v$ --- even if it would be more adequate, comparing with the Riemannian case, to call it the square of the norm.
An orthonormal frame is a frame
$(e_0, e_1, . . . , e_n)$
where $e_0$ as norm $−1$, every $e_i$ $(i \geq 2)$ has norm +1 and every scalar product $\langle e_i, e_j\rangle$ with $i \neq j$ is $0$.

We denote by $\SO_0(1,n)$, $\SO_0(2,n)$ the identity components of respectively
$\SO(1,n)$, $\SO(2,n)$ ($n \geq 2$). For any cocompact lattice $\Gamma$ of $\SO_{0}(1,n)$ and
any Lie group $G$ we denote by $\op{Rep}(\Gamma, G)$ the moduli space of representations\index{moduli space} of
$\Gamma$ into $G$ modulo conjugacy, equipped with the usual topology as an algebraic variety (see for example \cite{goldmillson}):
$$\op{Rep}(\Gamma, G) := \op{Hom}(\Gamma, G)/G.$$

Finally, if $(M_1, g_1)$ and $(M_2, g_2)$ are two pseudo-Riemannian manifolds, we denote by $g_1 - g_2$ the metric $p_1^*g_1 - p_2^*g_2$ on $M_1 \times M_2$, where $p_1$ and $p_2$ are the projections onto the first and second factor, respectively.
\medskip

\noindent\textbf{Index of notations.}
We introduce in this survey many objects and notions. We provide here for the reader's convenience an index of these objects, indicating the
page in which each of them is introduced (the list does not include objects already introduced such as $\SO_0(2,n)$ or $\cT_{2n}$).
\begin{itemize}
  \item $c$, $\hat{c}$, $\dot{c}$: causal curve, extension of the causal curve, derivative (beginning of Section \ref{sec:introcausal}).
  \item $[g]$: conformal class of the pseudo Riemannian metric $g$,
  \item $I^\pm$, $J^\pm$: future or past; causal future or past (Section \ref{sec:introcausal}).
  \item $U(p,q)$, $\overline{U}(p,q)$: diamond, closed diamond (Section \ref{sec:introcausal}).
  \item $P(S)$, $F(S)$, dev$(S)$: past development, future development, development of the closed edgeless achronal (CEA) subset $S$ (Section \ref{sec:introcausal}).
  \item $\mathcal C$: a category of spacetimes (Section \ref{sec:introcausal}).
  \item $M_{max}$: maximal extension of the spacetime $M$ (Theorem \ref{teo.clara}).
  \item $L(c)$, $d_{lor}$: length of the causal curve $c$, Lorentzian distance (just after Remark \ref{rk:GHconformal}).
  \item $(\RR^{r,s}, q_{r,s})$ :  $r+s$-dimensional vector space equipped with a quadratic form of signature $(r,s)$ (Section \ref{sec:models}).
  \item $\Ein^{1,n}$, $\overline{\Ein}^{1,n}$ and $\mathrm{p}: \wt\Ein^{1,n} \to \Ein^{1,n}$: Einstein universes of dimension $n+1$ and their universal covering (Section \ref{sub.einuniv}).
  \item $\delta$ and $\delta_0$: transformations on $\overline{\Ein}^{1,n}$ (Section \ref{sub.einuniv}).
  \item $E(\wt\Lambda)$, $E(\Lambda)$ : globally hyperbolic domains of $\wt\Ein^{1,n}$ or $\Ein^{1,n}$ associated with a CEA $\wt\Lambda$ or $\Lambda$ (Section \ref{sub.einuniv}).
  \item Fill$(\wt\Lambda)$ : filling of the CEA $\wt\Lambda$ (Section \ref{sub.einuniv}).
  \item $S(\mathcal{C}_{n+1})$ : Klein model of $\Ein^{1,n}$ (Section \ref{sub.einuniv}).
  \item Conv$(\Lambda)$ and Conv$^*(\Lambda)$ : convex hull of an achronal subset of $\Ein^{1,n}$ and its convex dual (Section \ref{sub.einuniv}).
  \item $I^+_0$ : future of the origin in Minkowski space (Section \ref{sub:mink}).
  \item Mink$^+(\tilde{x})$, $\mathcal I^\pm$ : Minkowski domain in $\wt\Ein^{1,n}$ associated with an element $\tilde{x}$ of $\wt\Ein^{1,n}$ and its Penrose components (Section \ref{sub:mink}).
  \item $\HH(\textsf x)$, $\partial\HH(\textsf x)$ : totally geodesic hypersurface of $\HH^{n+1}$ and its boundary defined by an element $\textsf x$ of de Sitter space $\dS^{1,n}$ (Section \ref{sub:dsbasic}).
  \item ${\mathcal B}(\SS^n)$ : space of round disks in $\SS^n$ (naturally identified with $\dS^{1,n}$, see Section \ref{sub:dsbasic}).
  \item $\DS^{1,n}$ : Klein model of de Sitter space (Section \ref{sub:dsbasic}).
  \item $\partial_{\pm}dS^{1,n}$ : past and future conformal boundaries of $\dS^{1,n}$ (Section \ref{sub:dsbasic}).
  \item $\ADS^{1,n}$ : Klein model of Anti-de Sitter space (Section \ref{sub:adsbasic}).
  \item $U(\textsf x)$ : affine domain centered at an element $\textsf x$ of $\ADS^{1,n}$ (Section \ref{sub:adsbasic}).
  \item $H^\pm(\textsf x)$ : past and future hyperplanes dual to an element $\textsf x$ of $\ADS^{1,n}$ (Section \ref{sub:adsbasic}).
  \item $\mathcal G$, $\mathcal K$ : Lie algebras of $\SO_0(2,n)$ and of its maximal subgroup $K$ (Remark \ref{rk:geodinvisible}).
  \item ${\mathcal U}^{1,n}$ : space of future oriented vectors tangent to $\AdS^{1,n}$ of norm $-1$ (Remark \ref{rk:geodinvisible}).
  \item $\lambda$ : Liouville form on ${\mathcal U}^{1,n}$ (Remark \ref{rk:geodinvisible}).
  \item $J$, $\omega$ : complex structure and K\"ahler form on $\cT_{2n}$ (Remark \ref{rk:geodinvisible}).
  \item $\nu$, $B$, II : Gauss map, shape operator and second fundamental form of a smooth spacelike surface in $\AdS^{1,n}$ (end of Section \ref{sub:adsbasic}).
  \item $M_\Lambda(\Gamma) = \Gamma\backslash\Omega(\Lambda)$ :  model maximal globally hyperbolic flat spacetime (Section \ref{sub:flatGH}).
  \item $M(\Sigma)$ : maximal globally hyperbolic de Sitter spacetime associated with the $(\SS^n, \SO_0(1,n+1))$-manifold $\Sigma$ (Section \ref{sub:dSGH}).
  \item ${\mathcal H}^\pm(\Lambda)$ : past and future horizons of a globally hyperbolic domain $\Omega(\Lambda)$ of $\AdS^{1,n}$ (Section \ref{sec.regads}).
  \item $E^-_0(\Lambda)$ : Past tight region of $\Omega(\Lambda)$ (Section \ref{sec.regads}).
  \item $E(\Lambda_{k, \ell})$ : Split AdS-spacetime (end of Section \ref{sec.regads}).
  \item $D(\Lambda)$ : conformal boundary of the invisible domain $\Omega(\Lambda)$ of $\AdS^{1,n}$ (Section \ref{sub:btz}).
  \item $\Lambda_\Gamma$ : limit set of an achronal group of isometries (Section \ref{sec:discrete}).
  \item eu$_b(\Gamma)$ : bounded Euler class of the subgroup $\Gamma$ of $\SO_0(2,n)$ (Theorem \ref{teo:euler}).
  \item $\op{Rep}_{an}(\Gamma, \SO_0(1,G))$ : space of Anosov representations of $\Gamma$ into $G$ (Section \ref{sub:anosov}).
  \item $\rho_L$, $\rho_R$ : left and right representations associated with a MGHC $\AdS^{1,2}$ spacetime (beginning of Section \ref{sub:ads3}).
  \item $\lambda^\pm$ : pleating laminations on the boundary of the convex core (Section \ref{sub:ads3}).
\end{itemize}


\section{A brief introduction to causality notions}
\label{sec:introcausal}
Let $(M^{n+1}, g)$ be a Lorentzian manifold. A tangent vector is \emph{spacelike}\index{spacelike vector} if its norm is positive; \emph{timelike}\index{timelike vector} if its norm
is negative; \emph{lightlike}\index{lightlike vector} if it is non zero and its norm is $0$. We also define \emph{causal\index{causal vector}} vectors as tangent vectors that are timelike or lightlike. An immersed hypersurface is \emph{spacelike}\index{spacelike hypersurface}
if all vectors tangent to $S$ are spacelike; it is \emph{nontimelike}\index{nontimelike hypersurface} if tangent vectors are all
spacelike or lightlike.
A \emph{causal}\index{causal curve} (resp. \emph{timelike}\index{timelike curve}) curve is an immersion $c : I \subset \mathbb{R} \to M$
such that for
every $t$ in $I$ the derivative $\dot{c}(t)$ is causal (resp. timelike). This notion extends
naturally to non-differentiable curves (see below, or \cite{beem}). Such a curve is \emph{extendible} if there
is another causal curve $\hat{c} : J \to M$ and a homeomorphism $\phi : I \to K \varsubsetneq J$ such
that $c$ coincides with $\hat{c} \circ \phi$. The causal curve c is \emph{inextendible}\index{inextendible causal curve} if it is not
extendible.
\medskip

\noindent\textbf{Conformal Lorentzian manifolds.}\index{conformal Lorentzian manifold}
The notion of timelike, lightlike and causal vectors or curves are the same for Lorentzian metrics in the same conformal class. Therefore, all the causality notions to be presented below apply to conformally Lorentzian manifolds $(M, [g])$, where $[g]$ denotes the conformal class of the Lorentzian metric $g$.
\medskip

\noindent\textbf{Time orientation.} We always assume that the manifold $M$ is oriented.
On $(M, [g])$ we have another orientability notion: a \emph{time orientation}\index{time orientation} of $(M, [g])$ is a continuous
choice, for every $p$ in $M$, of one of the two connected components of the set of timelike vectors at $p$.
When such a choice is possible, $(M, [g])$ is \emph{time-orientable,} and in short we will mention time orientable (conformal class of) Lorentzian manifolds as \emph{(conformal) spacetimes.}
Any conformal manifold is doubly covered by a time-orientable one.
Once the time-orientation has been selected
we have a notion of \emph{future-oriented}\index{future-oriented} or \emph{past-oriented}\index{past-oriented} causal vectors, therefore of
causal curves. We also have the notion of time function\index{time function}: a map $t: M \to \RR$ which is non decreasing along
any causal curve. Note that a time function may be non differentiable, and that a differentiable map
$f: M \to \RR$ is a time function if and only if its differential takes non negative values on future oriented
causal vectors.
\medskip

\noindent\textbf{Causality notions.} Two points in $M$ are \emph{causally related} if there exists a
causal curve joining them; they are \emph{strictly causally related} if moreover this
curve can be chosen timelike.
More generally: let $E$ a subset of $M$ and $U$ an open neighborhood of $E$ in $M$.
$E$ is \emph{achronal}\index{achronal subset} in $U$ if there is no timelike curve contained in
$U$ joining two points of the subset. It is \emph{acausal,}\index{acausal subset} or \emph{strictly achronal} in $U$ if there
is no causal curve contained in $U$ joining two points of $E$. We say simply that $E$
is (strictly) achronal if it is (strictly) achronal in $U = M$. Finally, we say that $E$
is locally (strictly) achronal if every point $p$ in $E$ admits a neighborhood $U$ in $M$
such that $E \cap U$ is (strictly) achronal in $U$. Spacelike hypersurfaces are locally acausal,
and nontimelike hypersurfaces are locally achronal.
\medskip

\noindent\textbf{Past, future.}
The \emph{future}\index{future} of a subset $A$ of $M$ is the open set $I^+(A)$ made of final points of future
oriented timelike curves not reduced to one point and starting from a point of $A$.
The \emph{causal future}\index{causal future} $J^+(A)$ of $A$ is the (non necessarily closed) set of final
points of future oriented causal curves,
possibly reduced to one point and starting from a point of $A$ (hence $A$ itself belongs
to its causal future). The (causal) past\index{past} \index{causal past} ($J^-(A)$) $I^-(A)$ of $A$ is the (causal) future of $A$ when the
time-orientation of $M$ is reversed.
This induces two partial orders on $M$: for every $p$ and $q$ in $M$, we write $p \preceq q$ if $q$
lies in the causal future of $p$, and $p \prec q$ if $q$ lies in $I^+(p)$.
\medskip

\noindent\textbf{Alexandrov topology.} An \emph{open diamond}\index{open diamond} is a domain $U(p, q) = I^-(p) \cap I^+(q)$
that is the intersection between the future and the past of two points $p$, $q$. Open diamonds
form the basis of some topology on M, the so-called
\emph{Alexandrov topology}\index{Alexandrov topology} (see \cite{beem}). Every $U(x, y)$ is open for the manifold
topology, but the converse in general is false; when it holds, $(M, [g])$ is said \emph{strongly causal.}\index{strongly causal}

Strong causality is equivalent to the following property (Proposition $3.11$ of \cite{beem}): for every point $p$ in $M$, every neighborhood of $p$ contains an open neighborhood $U$ (for the usual manifold
topology) of $p$ which is \emph{causally convex,}\index{causally convex} i.e. such that any causal curve in $M$ joining
two points in $U$ is actually contained in $U$.

From now on, we always assume that the spacetime $(M, [g])$ is strongly causal.

For any $p$, $q$ in $M$ the \emph{closed diamond}\index{closed diamond} $\overline{U}(p,q)$ is the intersection between the causal future of $p$ and the causal past of $q$.

\medskip

\noindent\textbf{Refined causality notions.} In strongly causal spacetimes, one can extend the class of causal curves in the following way:
a curve $c : I \subset \mathbb{R} \to M$ is causal and future oriented if it is locally non decreasing
for the usual order on $I$ and the partial order $\preceq$. It is strictly causal if it is locally increasing for
the partial order $\prec$. It is an easy exercise to see that acausal curves are locally Lipschitz, but
in general non smooth.

One also has the important following notion among locally achronal subsets: a locally achronal subset $A$ is \emph{edgeless}\index{edgeless achronal sunset} if every point $p \in A$ admits a neighborhood $U$ in $M$ such that any causal curve contained in $U$ and with extremities in respectively $I^+_U(A)$ and $I^-_U(A)$ crosses $A$. Then, closed edgeless achronal subsets (abbrev. CEA) are natural generalizations of (smooth) nontimelike hypersurfaces: they are locally graphs of Lipschitz maps.
\medskip

\noindent\textbf{Global hyperbolicity.}
A spacetime $(M, [g])$ is \emph{globally hyperbolic}\index{global hyperbolicity} (abbrev. GH) if:

–- it is strongly causal,

–- for any $p$, $q$ in $M$ the closed diamond $\overline{U}(p,q)$ is compact or empty.

This definition makes clear that, in a globally hyperbolic spacetime $(M, [g])$, an open domain $V$ of
$M$ is globally hyperbolic if and only if it is causally convex. Indeed,
each of these notions is equivalent to the fact that
for every $x$, $y$ in $V$, the closed diamond $\overline{U}(x,y)$ in $M$ coincides with the closed diamond in $V$.
It follows that intersections of GH domains of $(M, [g])$ are still GH.

The notion of global hyperbolicity is closely related to the notion of Cauchy
surfaces that we define now: let $S$ be a spacelike hypersurface embedded in $M$ (or, more generally, a CEA in $M$). The past development $P(S)$ (resp. the future development $F(S)$)
is the set of points $p$ in $M$ such that every inextendible causal path containing $p$
meets $S$ in its future (resp. in its past). The \emph{Cauchy development}\index{Cauchy development} dev$(S)$ is the union
$P(S) \cup F(S)$.
When dev$(S)$ is the entire $M$, $S$ is a \emph{Cauchy hypersurface.}\index{Cauchy hypersurface}
An important fact,
that can be considered as a generalisation of the Hopf-Rinow Theorem, is R. Geroch's Theorem (\cite{gerochdependence}):

\begin{teo}\label{teo:critereGH}
A strongly causal spacetime $(M, [g])$ is globally hyperbolic
if and only if it admits a Cauchy hypersurface. In this case, it is foliated by Cauchy hypersurfaces; more
precisely, there is a smooth time function $t: M \to \RR$ such that every level set of $t$ is a Cauchy hypersurface.
\end{teo}

It follows directly from this theorem that every GH spacetime $(M, g)$ is isometric
to a product $S \times \RR$ equipped with a metric of the form $\bar{g}_t - Ndt^2$ where
$\bar{g}_t$ is a one parameter family of Riemannian metrics on $S$ and $N: M \to ]0, +\infty[$ is a positive
function, called the \emph{lapse function} (see Proposition 6.6.8 of \cite{hawkingellis}).

In particular, Cauchy hypersurfaces in a given GH spacetime are diffeormophic to each other.

\begin{remark}
 There has been some imprecision in the literature concerning
the proof the smoothness of the splitting of globally hyperbolic spacetimes. See
\cite{sanchez, sanchez2, sanchezcausal} for a survey on this question and a complete proof of the smoothness
of the splitting $M \approx S \times \mathbb{R}$. See also a more recent proof with different methods, in \cite{fathi}.
\end{remark}

\begin{remark}
  The notion of global hyperbolicity in terms of Cauchy hypersurfaces has been introduced by J. Leray (1952).
  The key point is that the finite propagation property of the Einstein equations ensures that a metric
  solution of the Einstein equations is completely determined by its restriction to a neighborhood of a Cauchy hypersurface $S$ (more precisely, by the Riemannian metric obtained by restricting the Lorentzian metric to $S$, and
  by its second fundamental form). Therefore, globally hyperbolic spacetimes form a well posed problem from
  the viewpoint of Partial Differential Equations.
\end{remark}

\begin{remark}
  GH spacetimes are never compact. The suitable compactness notion is {spatial compactness}: a spacetime
  is \emph{globally hyperbolic spatially compact}\index{spatially compact} (abbrev. GHC) if it admits a compact Cauchy hypersurface --- all Cauchy hypersurfaces are then compact.
\end{remark}

\medskip

\noindent\textbf{Maximal globally hyperbolic spacetimes.}
An isometric embedding $f : M \to N$ is a \emph{Cauchy embedding} if
the image by $f$ of any Cauchy hypersurface of $M$ is a Cauchy hypersurface of N.

In this paragraph, we have to treat separately Lorentzian spacetimes and conformal Lorentzian spacetimes.
If $(M, g)$ and $(N, h)$ are GH Lorentzian spacetimes, a map $f : M \to N$ is a \emph{Cauchy embedding}\index{Cauchy embedding} if it is an isometric embedding such that the image by $f$ of any Cauchy hypersurface in $M$ is a Cauchy hypersurface in $N$. A \emph{conformal} Cauchy embedding is a conformal embedding $f : M \to N$ between conformal spacetimes
mapping Cauchy hypersurfaces into Cauchy hypersurfaces. Note that a conformal Cauchy embedding might be non-isometric, therefore not a
Cauchy embedding in our terminology.

Let $\mathcal C$ be a category of Lorentzian spacetimes, i.e. a class of Lorentzian spacetimes stable by isometries, by union, and restriction to open domains --- for example, the category of $C^r$ spacetimes, or the category of analytic spacetimes etc$\ldots$
A GH $\mathcal C$-spacetime $(M,g)$ is $\mathcal C$-\emph{maximal} (abbrev. $\mathcal C$-MGH)\index{maximal globally hyperbolic} if every Cauchy
embedding $f : M \to N$ in a $\mathcal C$-spacetime $N$ is surjective (hence a global isometry).

We have a similar notion of maximal conformal GH spacetimes, but where $\mathcal C$ is a category of
conformal spacetimes. The two notions may differ; as we will see later, a Lorentzian spacetime may be
maximal among spacetimes of constant curvature, but not maximal in the category of conformally flat
spacetimes.

A Lorentzian category $\mathcal C$ is \emph{rigid}\index{rigid category of spacetime} if it has the following property: if $p$, $q$ are points in $\mathcal C$-spacetimes $M$, $M'$ such that any isometry between $J^-(p)$ and $J^-(q)$ extends to an isometry between neighborhoods of $p$ and $q$. The traditional example is the category of solutions of the Einstein equations
in the void, or for us, the category of constant curvature spacetimes.

One has a similar notion of \emph{rigid conformal categories}: the ones for which any conformal
diffeomorphism between $J^-(p)$ and $J^-(q)$ extends to a conformal diffeomorphism between neighborhoods of $p$ and $q$.

\begin{teo}[\cite{clara}]\label{teo.clara}
Let $\mathcal C$ be a rigid category of Lorentzian spacetimes or of conformal spacetimes. Then any GH $\mathcal C$-spacetime $M$ admits a Cauchy embedding $f : M \to M_{max}$ in a $\mathcal C$-MGH spacetime. Moreover, $M_{max}$ is unique up to right composition by an isometry in the case of Lorentzian categories, and up to
right composition by a conformal diffeomorphism in the case of conformal categories.
\end{teo}

\begin{remark}\label{rk:GHconformal}
  This theorem was first established in the case of solutions of the Einstein equations (\cite{choquet}).
  There is a more recent proof, with the same ideas, but a different order allowing to avoid in a clever way the use of the Zorn lemma (\cite{sbierski}). In her work \cite{clara}, C. Rossi Salvemini observed that the proof applies tothe more general case of rigid categories of Lorentzian spacetimes, and also in the context of conformal Lorentzian categories that
  she introduced. Moreover, she proposed an entirely new proof, based on the notion of \emph{shadows} (intersections $J^\pm(p) \cap S$ between a Cauchy hypersurface $S$ and past/future of points). She also proved the following important result: \emph{if a GHC spacetime has nonpositive constant curvature and is maximal among spacetimes of constant curvature, then it is also maximal among conformally flat spacetimes.} This statement is false in the positive constant curvature case: a GH spacetime with constant positive curvature is \emph{never} maximal as a conformally flat spacetime.
\end{remark}

\noindent\textbf{Lorentzian distance.} Let $M$ be a time-oriented Lorentzian spacetime. The \emph{length-time} \index{length-time} $L(c)$ of a
causal curve $c : I \to M$ is the integral
over $I$ of the square root of $-\langle c′(t) | c′(t) \rangle$.
Observe that this is well-defined, since causal curves are always Lipschitz.
The \emph{Lorentzian distance}\index{Lorentzian distance} $d_{lor}(p, q)$ between two points $p$, $q$
is Sup$\{L(c) / c \in C(p, q)\}$ where $C(p, q)$ is the set of causal curves with extremities
$p$, $q$ (see for example \cite{cosmic}). By convention, if $p$, $q$ are not causally related, $d_{lor}(p, q) = 0$: when $M$ is globally hyperbolic, it defines
a continuous function $d_{lor}: M \times M \to [0, +\infty[$ since if $q$ lies on the boundary of $J^\pm(p)$ then there is a lightlike
curve joining $p$ to $q$ and $d_{lor}(p, q) = 0$.

\begin{theo}[Corollary $4.7$ and Theorem $6.1$ of \cite{beem}]
  If $M$ is globally hyperbolic, then $d_{lor} : M \times M \to [0,+\infty]$ is continuous and admits only finite values. Moreover, if $p$ is in the causal future of $q$, then there exists a geodesic $c$ with extremities $p$, $q$ such that $L(c) = d(x, y)$.
\end{theo}

It is to obtain this theorem that one does not restrict the definition of causal curves to piecewise $C^1$ curves.
\medskip

\noindent\textbf{Cosmological time.} In any spacetime, we can define the notion of cosmological time\index{cosmological time} (\cite{cosmic}):
For any $p$ in $M$, the cosmological time $\tau(p)$ is Sup$\{L(c) / c \in \mathcal R(p) \}$, where $\mathcal R(p)$ is the set of past-oriented causal curves starting at $p$.
This function could have in general a bad behavior: for example, in Minkowski
space, the cosmological time is everywhere infinite.

\begin{defi}\label{def:cosmic}
  A Lorentzian spacetime $(M,g)$ is said to have regular cosmological time\index{regular cosmological time} if:

–-- $M$ has finite existence time, i.e. $\tau(p) < \infty$ for every $p$ in $M$,

–-- for every past-oriented inextendible curve $c : [0,+\infty[\to M,$ we have lim$_{t\to\infty} \tau(c(t)) = 0$.
\end{defi}

Theorem $1.2$ in \cite{cosmic} expresses many nice properties of spacetimes with regular
cosmological time functions. We need only the following statement:

\begin{theo}\label{th:cosmic}
  If $M$ has regular cosmological time, then the cosmological time is
Lipschitz regular and $M$ is globally hyperbolic.
\end{theo}

\section{Model spacetimes}
\label{sec:models}
In this section, we describe the model spacetimes for every sign of the (constant) curvature. This includes
a description of their causal curves and their achronal subsets. We will end by the presentation of the space $\cT_{2n}$ of timelike geodesics in the anti-de Sitter space that will play an important role in this survey.
In the next section, we will use this
material for the classification of maximal GH spacetimes of constant curvature.

An important feature is that causality notions in the model spacetimes are much easier to deal with once it is observed that they all admit a conformal embedding in the Einstein Universe. We therefore start with this
central geometric model. As a reference for the content of this section, we mention \cite{francesthesis, theseclara} as very complete references (in French), the papers \cite{franceseinstein, clara} extracted from these works, and also \cite{merigot, ABBZ}.
\medskip

\noindent\textbf{On notations:} all model spacetimes  involve $(\RR^{r,s}, q_{r,s})$, for some integers $r$, $s$, where
$\RR^{r,s}$ is the $r+s$-dimensional vector space $\RR^{r,s}$ and $q_{r,s}$ a quadratic form of signature $(r,s)$. More precisely, the first $r$ coordinates of $\RR^{r,s}$ are denoted by
$u_1, \ldots, u_r$ and the other coordinates by $x_1, \ldots, x_s$. Elements of $\RR^{r,s}$ are denoted by $\textsf{x}$, $\textsf y$, $\ldots$ The quadratic form is then:
$$q_{r,s}(\textsf x) = -u_1^2 -\ldots - u_r^2 + x_1^2 + \ldots + x_s^2$$
The associated scalar product is denoted by $\langle . \mid . \rangle_{r,s}$, or simply $\langle . \mid . \rangle$.

We will also denote by $(\SS^n, \bar{g}_n)$ the sphere of dimension $n$ equipped with its usual metric
$\bar{g}_n$: the restriction of $q_{0,n+1}$ to the unit sphere $\{ q_{0,n+1}=1 \}$ of $\RR^{0,n+1}$.
The distance on $\SS^n$ induced by $\bar{g}_n$ is denoted by $d_n$.

\subsection{Einstein universe}
\label{sub.einuniv}
The four-dimensional Einstein universe was the first cosmological model for our universe proposed by A. Einstein soon after the birth of General Relativity. The $n+1$ dimensional Einstein universe $\wt\Ein^{1,n}$ \index{Einstein universe} can be simply described as the (oriented) product $\SS^n \times \RR$ of the $n$-dimensional sphere and the real line, equipped with the metric $\bar{g}_n - dt^2$, and time-oriented so that the coordinate $t$
is a time function.

The importance of the Einstein universe is essentially due to the following extension of the Liouville Theorem (\cite{francesliouville}): \emph{when $n\geq 2$, any conformal transformation between two open subsets of $\wt\Ein^{1,n}$
extends to a global conformal transformation}. It follows that conformal Lorentzian spacetimes of dimension $\geq 2+1$ are locally modeled on the Einstein universe. Hence we don't really consider the Einstein universe as a Lorentzian manifold, but as a conformally Lorentzian spacetime.

We also consider the product $\Ein^{1, n} = \SS^n \times \SS^1$ equipped with the $\bar{g}_n - \bar{g}_1$.
We denote by $\mathrm{p}: \wt\Ein^{1,n} \to \Ein^{1,n}$ the
cyclic covering map. Let $\delta: \wt\Ein^{1,n} \to \wt\Ein^{1,n}$ be the map $(\bar{x,} t) \mapsto (\bar{x,} t + 2\pi)$: it generates the Galois group of $\mathrm{p}$. We will also
consider the quotient $\overline{\Ein}^{1,n}$ of $\wt\Ein^{1,n}$ by $\delta_0: (\bar{x,} t) \mapsto (-\bar{x,} t+\pi)$,
even if its topology is slightly more difficult to handle. This quotient, which is doubly covered by
 $\Ein^{1,n}$, is sometimes called \emph{Einstein universe} in the literature.\medskip

\noindent\textbf{Photons.} Lightlike geodesics of $\Ein^{1,n}$, $\wt\Ein^{1,n}$ or $\overline\Ein^{1,n}$, considered
as non-parameterized curves, do not depend on the choice of the representative of the conformal class.
They are called \textit{photons.} \index{photon} In $\wt\Ein^{1,n}$, they are represented by the curves $t \mapsto (\bar{x}(t), t)$, where $t \mapsto \bar{x}(t)$ is a geodesic of $(\SS^n, \bar{g}_n)$. For every $\tilde{x}$ in $\wt\Ein_n,$ $\delta_0(x)$ is its first conjugate point in the future: by this, we mean that
every future oriented photon exiting from $\tilde{x}$ also contain $\delta_0(\tilde{x})$, and that there is no other such conjugate point in $J^+(\tilde{x}) \cap J^-(\delta_0(\tilde{x}))$.

The union of all photons containing $\tilde{x}$ is the \emph{lightcone}\index{lightcone} $C(\tilde{x})$. If we write $\tilde{x}$ as a pair $(\bar{x}, t)$ in $\SS^n \times \RR$, the lightcone $C(\tilde{x})$ is the set
of pairs $(\bar{y,}s)$ such that the difference $t-s$ is equal to $d_n(\bar{x,}\bar{y})$ modulo $2\pi$.
It is a cylinder pinched at every conjugate points $\delta_0^k(\tilde{x})$ (see for example Figure $4.3$ in \cite{francesthesis}).

The projections of lightcones in $\Ein^{1,n}$ are also called lightcones; a lightcone $C(x)$ for $x$ in
$\Ein^{1,n}$ has two singular points, $x$ and its conjugate $\delta_0(x)$.
\medskip

\noindent\textbf{Causal curves. }More generally, causal curves in $\wt\Ein^{1,n}$, suitably parameterized, are the curves $t \mapsto (\bar{x}(t), t)$, where $t \mapsto \bar{x}(t)$ is a $1$-Lipschitz map from an interval $I$ into $\SS^n$ --- it is timelike if $\bar{x}$ is contracting, i.e.:
$$\forall s,t \in \RR \;\; d_n(\bar{x}(s), \bar{x}(t)) < \mid s-t \mid.$$
In particular, inextendible causal curves are the ones parameterized by $I = \RR$. It clearly follows that $\wt\Ein^{1,n}$ is globally hyperbolic spatially compact: every level set $\{ t = \rm Const. \}$ is a Cauchy hypersurface.
\medskip

\noindent\textbf{Achronal subsets.}
Achronal subsets of $\wt\Ein^{1,n}$ are precisely the graphs
of $1$-Lipschitz functions $f: \Lambda_0 \rightarrow {\mathbb R}$ where
$\Lambda_0$ is a subset of ${\mathbb S}^{n}$. The achronal set is acausal if and only
if $f$ is $1$-contracting. It is closed if and only if $\Lambda_0$ is closed, and edgeless if and only if
$\Lambda_0$ is open.
In particular, closed edgeless achronal subsets are
exactly the graphs of the $1$-Lipschitz functions $f:\SS^{n}\to\RR$: they are
topological $n$-spheres, which are all Cauchy hypersurfaces.

\emph{Stricto sensu,\/} there is no achronal subset in $\Ein^{1,n}$ since closed timelike curves
through a given point cover the entire
$\Ein^{1,n}$. Nevertheless, we can keep track of this notion in $\Ein^{1,n}$ by defining ``achronal'' subsets
of $\Ein^{1,n}$ as projections of genuine achronal subsets of $\wt\Ein^{1,n}$. This definition is
justified by the following results (Lemma 2.4, Corollary 2.5 in \cite{merigot}):
\emph{The restriction of $\mathrm{p}$ to any achro\-nal subset of \/$\/\wt\Ein^{1,n}$ is injective.}
\emph{Moreover, if  \/$\wt\Lambda_{1}$, $\wt\Lambda_{2}$ are two achronal subsets of \/$\wt\Ein^{1,n}$
admitting the same projection in \/$\Ein^{1,n}$, then there is an integer $k$ such that:}
$$\wt\Lambda_{1}=\delta^{k}\wt\Lambda_{2}$$
In this setting, closed edgeless achronal subsets of $\Ein^{1,n}$ are graphs of $1$-Lipschitz maps from
$\SS^n$ into $\SS^1$.
\medskip

\noindent\textbf{Globally hyperbolic domains.} Let $\wt\Lambda$ be a closed achronal subset of $\wt\Ein^{1,n},$ \ie the graph of a $1$-Lipschitz map $f: \Lambda_0 \to \RR$ where $\Lambda_0$ is a closed subset of $\SS^n$. Define two functions
$f^{-}, f^{+}:{\SS}^{n}\to\RR$ as follows:
\begin{align*}
f^{+}(\bar{x}) &:= \mbox{Inf}_{\bar{y} \in \Lambda_0} \{ f(\bar{y})+d_n(\bar{x,}\bar{y}) \} ,\\
f^{-}(\bar{x}) &:= \mbox{Sup}_{\bar{y} \in \Lambda_0} \{ f(\bar{y})-d_n(\bar{x,}\bar{y}) \} ,
\end{align*}
$f^+$ (respectively $f^-$) is the maximal (respectively minimal) $1$-Lipschitz map from $\SS^n$ into $\RR$ that coincides with $f$ on $\Lambda_0$.

Then, the set of points of $\wt\Ein^{1,n}$ which are not causally related to any point of $\wt\Lambda$ is:
$$
E(\wt\Lambda)=\{(\bar{x,} t)\in {\SS}^{n} \times \RR \mid
f^{-}(\bar{x})< t <f^{+}(\bar{x})\}.
$$

Observe that if $\bar{x}$ and $\bar{y}$ are two points in $\Lambda_0$ such that $|f(\bar{x}) - f(\bar{y})|= d_n(\bar{x}, \bar{y})$, then
the restrictions of $f^+$ and $f^-$ to any minimizing $d_n$-geodesic segment between $\bar{x}$ and $\bar{y}$ coincide. Let Fill$(\Lambda_0)$ be
the union of $\Lambda_0$ with the union of all minimizing $d_n$-geodesic segments joining two elements $\bar{x}$, $\bar{y}$ of $\Lambda_0$ such that $|f(\bar{x}) - f(\bar{y})|= d_n(\bar{x}, \bar{y})$, and let Fill$(f)$ be the restriction of $f^\pm$ to Fill$(\Lambda_0)$: the graph of Fill$(f)$ is a CEA of $\wt\Ein^{1,n}$ that we denote by Fill$(\wt\Lambda)$ and that we call the \textit{filling of $\wt\Lambda$}\index{filling of an achronal subset}. Then we have $E(\wt\Lambda) = E($Fill$(\wt\Lambda))$. In other words, we can restrict ourselves to \emph{filled CEA,} \ie CEA equal to their own fillings (see \cite[Remark 3.19]{quasifuchs}).

The domain $E(\wt\Lambda)$ may be empty, but exactly in the case where the filling Fill$(\wt\Lambda)$ is the entire sphere $\SS^n$.
A particular case when this happens is the case where $\wt\Lambda$ is \emph{purely lightlike,} i.e. the
case where $\Lambda_0$ contains two antipodal points $\bar{x}_0$ and $-\bar{x}_0$
such that the equality $f(\bar{x}_0)=f(-\bar{x}_0) +\pi$ holds (Lemma 3.6 in \cite{merigot}).
Purely lightlike achronal subsets are precisely the ones admitting as filling the union of lightlike geodesics joining two antipodal points of $\wt\Ein_n$.

If non-empty, the invisible domain $E(\wt\Lambda)$ is globally hyperbolic (indeed, it is easy to see that for any $p$, $q$ in $E(\wt\Lambda)$, the closed diamond $\overline{U}(p,q)$ is contained in $E(\wt\Lambda)$). More precisely,
the Cauchy hypersurfaces of $E(\wt\Lambda)$ are precisely
graphs $\Lambda_F$ of $1$-Lipschitz maps $F: \SS^n \setminus$ Fill$(\Lambda_0) \to \RR$ such that the extension of $F$ to
Fill$(\Lambda_0)$ coincides with Fill$(f)$.

In the limit case $\Lambda_0 = \emptyset$, we have $f^+ = +\infty$ and $f^- = -\infty$, and the Cauchy development of the graph of any $1$-Lipschitz map $F: \SS^n \to \RR$ defined on the entire sphere is the entire Einstein universe $\wt\Ein^{1,n}$.

In summary, the theory of globally hyperbolic domains of $\wt\Ein^{1,n}$ coincides with the theory
of $1$-Lipschitz maps on $\SS^n$.

As for achronal subsets, even if $\Ein^{1,n}$ is not strongly causal, one can abusively project
the notion of globally hyperbolic domains into $\Ein^{1,n}$, thanks to the following lemma:

\begin{lem}
\label{le.one-to-one}
For every (non-empty) closed achronal set
$\wt\Lambda\subset\wt\Ein^{1,n}$, the projection of $E(\wt\Lambda)$ onto $E(\Lambda) = \op{p}(E(\wt\Lambda))$ is one-to-one.
\end{lem}
\medskip

\noindent\textbf{Klein model.} Einstein universe can also be defined in the following way: let $\mathcal C_{n+1}$ be the null-cone in $(\RR^{2,n+1,} q_{2,n+1})$, and let
$\SS(\mathcal C_{n+1})$ be its projection in the space $\SS(\RR^{2,n+1})$ of rays in $\RR^{2,n+1}$.
$\SS(\RR^{2,n+1})$ is a double covering of the usual projective space $\mathbb P(\RR^{2,n+1})$, therefore we call it,
slightly abusively, the projectivization of $\RR^{2,n+1}$. Observe that the convexity is well defined
in $\SS(\RR^{2,n+1})$: one can simply define convex subsets of $\SS(\RR^{2,n+1})$ as radial projections
of convex cones in $\RR^{2,n+1}$. In particular, convex hulls $\op{Conv}(B)$ of subsets $B$ of $\SS(\RR^{2,n+1})$,
in particular, of subsets of $\SS(\mathcal C_{n+1})$,
are well defined (but we don't mean that they are contained in $\SS(\mathcal C_{n+1})$!).

The quadratic form $q_{2,n+1}$ induces a natural conformally Lorentzian structure on $\SS(\mathcal C_{n+1})$. More precisely, for any section $\sigma: \SS(\mathcal C_{n+1}) \to \mathcal C_{n+1}$, the pull-back $\sigma^\ast q_{2,n+2}$ is a Lorentzian metric $g_\sigma$, and the conformal class $[g_\sigma]$ does not depend on $\sigma$. This conformally Lorentzian metric happens to be  conformally isometric
to $(\Ein^{1,n}, \bar{g}_n - \bar{g}_1)$.

The pair $(\SS(\mathcal C_{n+1}), [g_\sigma])$ is the \emph{Klein model}\index{Klein model of Einstein universe} of $\Ein^{1,n}$.
\medskip

\noindent\textbf{Isometry group.}
It is clear from the Klein model that the group of conformal transformations of $\Ein^{1,n}$ preserving the orientation and the time orientation is $\SO_0(2, n+1)$. Let $\wt{\SO}_0(2,n+1)$ be the group of conformal transformations of $\wt\Ein^{1,n}$ preserving the orientation and the time orientation. There is a natural projection from $\wt{\SO}_0(2,n+1)$ into $\SO_0(2, n+1)$ whose kernel is spanned by the transformation $\delta$ generating the Galois group of $\op{p}: \wt\Ein^{1,n} \to \Ein^{1,n}$ defined previously. Therefore, there is a central exact sequence:
$$1 \to \ZZ \to \wt{\SO}_0(2,n+1) \to \SO_0(2,n+1) \to 1$$
Observe that $\widetilde{\SO}_0(2,n+1)$ is not the universal covering of $\SO_0(2,n+1)$: there is a retraction of $\SO_0(2, n+1)$ onto
$\SO(2) \times \SO(n)$ hence the fundamental group of $\SO(2,n+1)$ is not cyclic but isomorphic\footnote{We thank
the referee to have pointed out this fact to us.} to $\ZZ \times (\ZZ/2\ZZ)$.

\begin{remark}
Concerning the notation: in the sequel, we always have in mind the identifications
$\Ein^{1,n} \approx \SS(\cC_{n+1})$, and we frequently switch from one model to the other.
We denote by $x$ elements of $\Ein$, using the notation $\op{x}$ when we want to insist
on the Klein model.
\end{remark}
\medskip

\noindent\textbf{Causality notions in the Klein model.}
Two elements $\op{x}$, $\op{y}$ of $\Ein$ are causally related if and only if
$\langle \op{x} \mid \op{y} \rangle \geq 0$. In particular, a subset $\Lambda \subseteq \Ein$ is achronal (respectively acausal) if and only if for every distinct $\op{x,} \op{y}\in\Lambda$ the
scalar product $\langle \op{x} \mid \op{y} \rangle$ is non-positive
(respectively negative).

Photons are projections on $\SS(P)$ of isotropic $2$-planes of $\RR^{2, n+1}$. The lightcone $C(\op{x})$ of
a point $\op{x}$ is the projection of $\mathcal C_{n+1} \cap \textsf{x}^\perp$, where $\textsf{x}^\perp$ is
the $q_{2,n+1}$-orthogonal of any representative $\textsf x$ of $\op{x}$.

Finally, for every achronal subset $\Lambda$ of $\Ein^{1,n} \approx \SS(\cC_{n+1})$,
the invisible domain $E(\Lambda)$ is:
$$E(\Lambda) = \{ \op{x} \in \SS(\mathcal C_{n+1}) \;\mid\; \forall \op{y} \in \Lambda \;\; \langle \op{x} \mid \op{y} \rangle < 0\}.$$
Recall that the dual of a convex subset $\SS(C)$ of $\SS(\RR^{2,n})$ is:
$$C^* = \{ \op{x} \in \SS(\mathcal C_{n+1}) \;\mid\; \forall \op{y} \in C \;\; \langle \op{x} \mid \op{y} \rangle \leq 0\}.$$
Hence, $E(\Lambda)$ coincides with the interior of $\SS(\mathcal C_{n+1}) \cap \op{Conv}^\ast(\Lambda)$, where $\op{Conv}^\ast(\Lambda)$ is the convex subset of $\SS(\RR^{2,n})$ dual to the convex hull $\op{Conv}(\Lambda)$: it is the intersection between a quadric and a convex subset of the projective space $\SS(\RR^{2,n})$ (cf. \cite[section $3.3$]{quasifuchs}).

\subsection{Minkowski space}\label{sub:mink}
For this section, we indicate as references \cite{barflat, bonsante, bardomaine}.
The \noindent\textbf{Minkowski space}\index{Minkowski space} is the affine space of dimension $n+1$ equipped
with the quadratic form $q_{1,n}$ on the underlying vector space $\mathbb{R}^{1,n}$ (for $n\geq 1$).
We slightly abuse notations, denoting it by $\mathbb{R}^{1,n}$, whereas it should really considered as an affine space, and not linear. We also use a coordinate system $(t, x_1, \ldots, x_n)$ such that:
$$q_{1,n}(\textsf x) = -t^2 + x_1^2 + \ldots + x_n^2.$$

The causal structure of the Minkowski space is particularly simple, because of its affine structure. It is convenient
to see it as the product $(\RR \times \RR^n, -dt^2 + |d\bar{\textsf x}|^2)$ of the line $\RR$ ``of time'' and the
Euclidian plane $\RR^n$, whose elements are denoted by $\bar{\textsf x}$. A time orientation is obtained by requiring the time coordinate $t$ to be a time function.

Let us fix an origin $0$, identifying Minkowski space
with its underlying linear space $\RR^{1,n}$. Let $I^+_0$ be the set of future-oriented
timelike tangent vectors at $0$. Then $I^+_0$ coincides with the future $I^+(0)$ through the canonical identification between $T_0\RR^{1,n}$ and $\RR^{1,n}$. Two elements $\textsf{x}$, $\textsf y$ are causally related
if and only if $\textsf y - \textsf x$ is timelike, and $\textsf y \in I^+(\textsf x)$ if and only if
$\textsf y - \textsf x$ lies in $I^+_0$.

The isometry group (as always, preserving all orientations) is the Poincar\'e group Isom$_0(\RR^{1,n})$, isomorphic
to $\SO_0(1,n) \ltimes \RR^{1,n}$.

Concerning the causality notions, almost all the discussion above in the case of Einstein universe applies, replacing
the sphere $(\SS^n, \bar{g}_n)$ by the Euclidean plane $(\RR^{0,n}, q_{0,n})$: up to reparametrization, causal (resp. timelike) curves  are maps $t \in I \subseteq \RR \mapsto \bar{\textsf{x}}(t)$ where
$\bar{\textsf{x}}: I \to \RR^{0,n}$ is
$1$-Lipschitz (resp. $1$-contracting). For inextendible curves we have $I = \RR$. Every horizontal hyperplane $\{ t = \rm Const. \}$ is a Cauchy surface. The geodesics of Minkowski space are affine lines. The achronal (resp. acausal) subsets are graphs $\Lambda_f$ of $1$-Lipschitz (respectively \/$1$-contracting) maps $f: \Lambda_0 \subseteq \RR^{0,n} \to \RR$, closed edgeless achronal subsets are graphs $\Lambda_f$ of $1$-Lipschitz maps defined
on the entire $\Lambda_0 = \RR^{0,n}$ --- but as we will se later, it is not always a Cauchy hypersurface for $\RR^{1,n}$.
\medskip

\noindent\textbf{Conformal model.} Minkowski space admits a conformal embedding in Einstein universe; actually,
it is conformally isometric to the complement in $\overline\Ein^{1,n}$ of any lightcone. It is actually more convenient to see it as one connected component of the complement of the lightcone $C(\tilde x_0)$ of some point $\tilde x_0$ in $\wt\Ein^{1,n}$: the ``extreme diamond'', intersection $I^+(\delta_0^{-1}(\tilde x_0)) \cap I^-(\delta_0^{-1}(\tilde x_0))$, which we denote by Mink$^+(\tilde x_0)$. If can take $\tilde x_0 = (\bar x_0, 0)$, then
Mink$^+(\tilde x_0)$ is the set of points $(\bar x, t)$ such that $| t | < d_n(\tilde x, \tilde x_0)$.

The boundary of Mink$^+(x_0)$ in $\Ein^{1,n}$ can therefore be seen as a \noindent\textbf{conformal boundary} of the Minkowski space, that has already been introduced by R. Penrose (\cite{Penrose}). It decomposes in several parts:

-- the point $\delta^{-1}_0(\tilde x_0)$, denoted by $i^-$ by Penrose: Mink$^+(\tilde x_0)$ is entirely contained in the future of $i^-$,

-- the point $\delta_0(\tilde x_0)$, also denoted by $i^+$: Mink$^+(\tilde x_0)$ is entirely contained in the past of $i^+$,

-- the point $i_0 := \tilde x_0$, called the ``spatial infinity'',

-- the complement of these three points is the union of two lightlike cylinders $\SS^{n-1} \times \RR$,
one in the future of Mink$^+(\tilde x_0)$ and denoted by $\mathcal{I}^+$, and the other, $\mathcal I^-$, in the past of
Mink$^+(\tilde x_0)$.

\begin{figure}[ht]
  \centering
  \includegraphics[width=8cm]{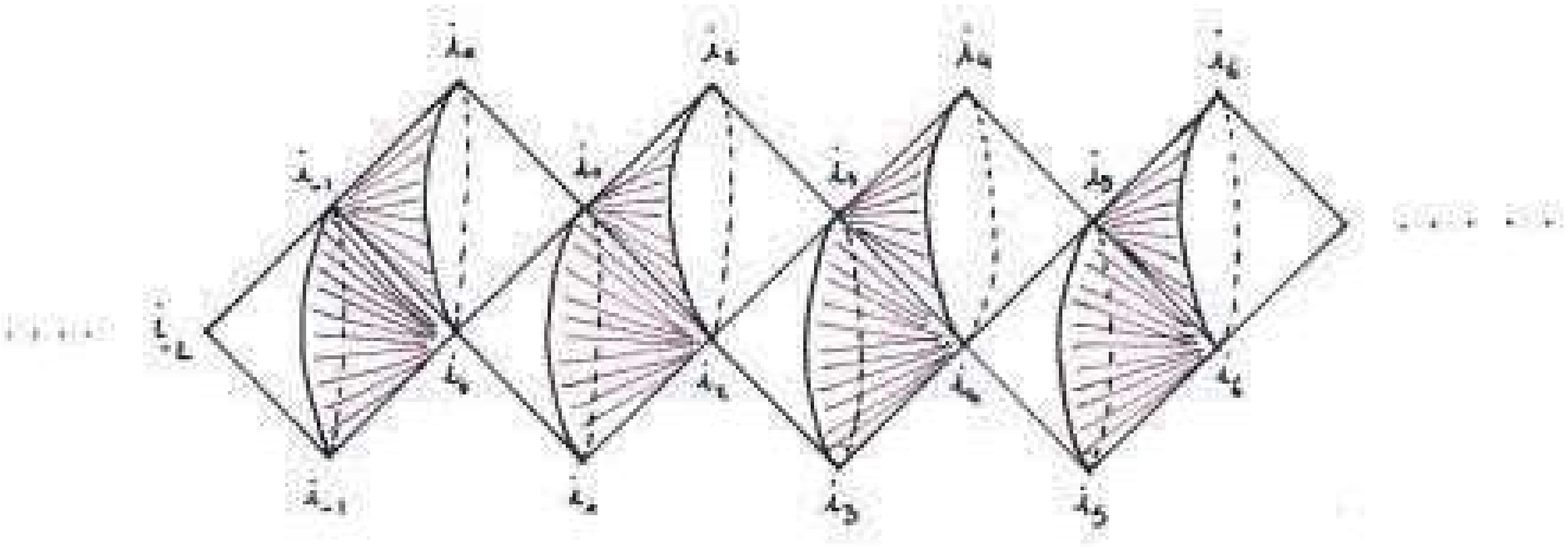}
\centering  \caption{Several Minkowski components.}
  \label{fig:minks}
\end{figure}

See Figure \ref{fig:minks} (borrowed from \cite[Figure $4.3$]{francesthesis}) where several Minkowski components Mink$^+(i_k)$
are depicted, where every $i_k$ is a point conjugate to $i_0 = \tilde x_0$ (\ie iterates of $\tilde x_0$ under $\delta_0$).

The intersection between Mink$^+(\tilde x_0)$ and the lightcone $C(\tilde x)$ of any point in $\mathcal I^\pm$ is
an affine hyperplane $H(\tilde x)$ in Mink$^+(\tilde x_0) \approx \RR^{1,n}$. More precisely: if $\tilde x \in \mathcal I^+$, then Mink$^+(\tilde x_0) \cap I^-(\tilde x)$ is the past of $H(\tilde x)$ in $\RR^{1,n}$, and if $\tilde x \in \mathcal I^-$, then
Mink$^+(\tilde x_0) \cap I^+(\tilde x)$ is the future of $H(\tilde x)$ in $\RR^{1,n}$. Therefore, $\mathcal I^-$ can be seen as the space of half affine Minkowski spaces, that are equal to their own future, and bounded by a lightlike hyperplane.

\subsection{De Sitter space}
\label{sub:dsbasic}
The \noindent\textbf{de Sitter space}\index{de Sitter space} $\dS^{1,n}$ is the hypersurface $\{ \textsf{x} \in
\RR^{1,n+1} /  \mathrm{q}_{1,n+1}(\textsf{x})=+1 \}$ endowed with the Lorentzian
metric obtained by restriction of $\mathrm{q}_{1,n+1}$. Hence, for the coordinates
$(t, x_1, \ldots , x_n)$ we have:
$$\dS^{1,n} := \{ (t, x_1, \ldots , x_n) \;\mid\; -t^2 + x_1^2 + \ldots + x_n^2 = +1\}.$$
We equip $\dS^{1,n}$ with the time orientation for which $t$ (or, more generally, for
every $\textsf u \in \HH^{n+1}$, the map $\textsf x \to -\langle \textsf x \mid \textsf u\rangle$) is a time function. The geodesics are the intersections between $\dS^{1,n}$ and $2$-planes in $\RR^{1, n+1}$.
Two points $\textsf{x}$, $\textsf{y}$ in $\dS^{1,n}$ are causally related if and only if
$\langle \textsf x \mid \textsf y\rangle \geq 1$.
\medskip

\noindent\textbf{Isometry group.}
The group of orientation preserving, time-orientation preserving isometries of $\dS^{1,n}$ is
$\SO_0(1,n+1)$.
\medskip

\noindent\textbf{Duality with the hyperbolic space.}
For every $\textsf x$ in $\dS^{1,n}$, the intersection $\textsf x^\perp \cap \HH^n$ is a
totally geodesic hypersurface $\HH(\textsf x)$ of hyperbolic space $\HH^{n+1}$. More precisely, $\{ \textsf y \in \HH^{n+1} \;\mid\; \langle \textsf y \mid \textsf x\rangle > 0\}$ is a half-hyperbolic space bounded by $\HH(\textsf x)$; in other words, $\dS^{1,n}$ can be seen as the space of \noindent\textbf{transversely oriented} half hyperbolic spaces.
\medskip

\noindent\textbf{De Sitter space is the space of round disks in $\SS^n$.}
Every $\HH(\textsf x)$ is characterized by its conformal boundary $\partial\HH(\textsf x)$ in
$\partial\HH^{n+1} \approx \SS^n$. Therefore, and this will be important later, $\dS^{1,n}$ can be identified with the space $\mathcal B(\SS^n)$ of round disks in $\SS^n$.

\medskip

\noindent\textbf{Klein model.}
The Klein model $\DS^{1,n}$ \index{Klein model of de Sitter space} is the projection of $\dS^{1,n}$ to $\SS(\RR^{1,n+1})$ \ie
$$\DS^{1,n} := \{ \op{x} \in \SS(\RR^{1,n+1}) \;/\; \langle \op{x} \mid \op{x} \rangle > 0 \}.$$

This is the complement in $\SS(\RR^{1,n+1})$ of the closures of two Klein models $\HH_\pm^{n+1}$:
$$\HH_\pm^{n+1} := \{ \op{x} \in \SS(\RR^{1,n+1}) \;/\; \langle \op{x} \mid \op{x} \rangle < 0, \;\; \pm t > 0 \}.$$
\medskip

\noindent\textbf{Conformal model }
De Sitter space $\dS^{1,n}$ is conformally equivalent to the domain
$dS^{1,n} := \SS^{n} \times ]-\pi, +\pi[ \; \subseteq \Ein^{1,n}$. It immediately follows
that $dS^{1,n}$ is a globally hyperbolic domain of $\Ein^{1,n}$, hence globally hyperbolic.
The achronal subsets are the graphs of $1$-Lipschitz maps $f: \Lambda_0 \subseteq \SS^n \to ]-\pi, +\pi[$.

The boundary of $dS^{1,n}$ has two components: the component $\SS^n \times \{ - \pi \}$ is the \emph{past conformal boundary} $\partial_-dS^{1,n}$, and the component $\SS^n \times \{ +\pi \}$ is the \emph{future conformal boundary} $\partial_+dS^{1,n}$. For every $x$ in $dS^{1,n}$, the intersection between the future (resp. the past) of $x$ in $\Ein^{1,n}$ and $\partial_+dS^{1,n}$ (resp. $\partial_-dS^{1,n}$) is a round ball $\mathcal B_+(\textsf x)$ (resp. $\mathcal B_-(\textsf x)$).
This is another way --- actually, \noindent\textbf{two} other ways --- to identify $\dS^{1,n}$ with $\mathcal B(\SS^n)$. For each of them, the causality relation in $dS^{1,n}$ corresponds to the inclusion relation in $\mathcal B(\SS^n)$.

\subsection{Anti-de Sitter space}
\label{sub:adsbasic}
\noindent\textbf{Anti-de Sitter space} \index{anti-de Sitter space} $\AdS^{1,n}$ is the hypersurface $\{ \textsf{x} \in
\RR^{2,n} /  \mathrm{q}_{2,n}(\textsf{x})=-1 \}$ endowed with the Lorentzian
metric obtained by restriction of $\mathrm{q}_{2,n}$.
We use a coordinate system $(u, v, x_1, \ldots , x_n)$ such that:
$$q_{2,n}(\textsf x) := -u^2 -v^2 + x_1^2 + \ldots + x_n^2.$$

We will also consider the coordinates $(r, \theta, x_{1,} \hdots, x_{n})$ with:
\[ u=r\cos(\theta), v=r\sin(\theta).\]
We equip $\AdS^{1,n}$ with the time orientation defined
by the vector field $\frac{\partial}{\partial\theta}$, \ie the time orientation such that the timelike vector field $\frac\partial{\partial\theta}$ is everywhere
future oriented.

Observe the analogy with the definition of Hyperbolic space $\HH^{n}$. Moreover, for every real number $\theta_0$, the subset $H_{\theta_0} := \{ (r, \theta, x_1,\hdots, x_n) /
\theta=\theta_{0} \} \subset \RR^{2,n}$ is a totally geodesic copy of  $\HH^{n}$ embedded in $\AdS^{1,n}$. More generally, the
totally geodesic subspaces of dimension $k$ in $\AdS^{1,n}$ are the
connected components of the intersections of $\AdS^{1,n}$ with the
linear subspaces of dimension $(k+1)$ in $\RR^{2,n}$.
In particular, geodesics are intersections with $2$-planes.
\medskip

\noindent\textbf{Conformal model }
Anti-de Sitter space $\AdS^{1,n}$ is conformally equivalent to the domain
$\DD^{n} \times \SS^1 \subseteq \Ein^{1,n}$, where $\DD^{n}$ is the open upper
hemisphere of $\SS^{n}$. The boundary of this domain is $\partial\DD^n \times \SS^1 \approx \SS^{n-1} \times \SS^1$,
hence conformally isometric to the Einstein universe $\Ein^{1,n-1}$ with one dimension less.
In other words, $\AdS^{1,n}$ is one of the two connected components of $\Ein^{1,n} \setminus \Ein^{1,n-1}$
for the natural conformal embedding of $\Ein^{1,n-1}$ in $\Ein^{1,n}$. In other words, one can
see $\Ein^{1,n-1}$ as the conformal boundary $\partial\AdS^{1,n}$.

As $\Ein^{1,n}$, $\AdS^{1,n}$ is not strongly causal; it contains many timelike geodesic loops.
But, its universal covering $\wt\AdS^{1,n}$, conformally
equivalent to $\DD^n \times \RR,$ is strongly causal. Its conformal boundary is
$\wt\Ein^{1,n-1}$.

However, $\wt\AdS^{1,n}$ is not globally hyperbolic: for example, diamonds $J^-(\bar x_0, t) \cap J^+(\bar{x}_0, 0)$ are not compact as soon as $t \geq 2\pi$.
\medskip

\noindent\textbf{Klein model.}
The Klein model $\ADS^{1,n}$ \index{Klein model of anti-de Sitter space} is the projection of $\AdS^{1,n}$ to $\SS(\RR^{2,n})$ \ie:
$$\ADS^{1,n} := \{ \op{x} \in \SS(\RR^{2,n}) \;/\; \langle \op{x} \mid \op{x} \rangle < 0 \}.$$

The topological boundary of $\ADS^{1,n}$ in $\SS(\RR^{2,n})$ is the Klein model $\SS(\mathcal C_{n})$
of $\Ein^{1,n-1}$. Observe that for any subset $\Lambda$ of $\Ein^{1,n-1} \approx \SS(\mathcal C_{n})$,
the convex hull $\op{Conv}(\Lambda)$ is contained in the conformal compactification $\AdS^{1,n} \cup \partial\AdS^{1,n}$ if and only if
$\Lambda$ is achronal in $\Ein^{1,n-1}$.
\medskip

\noindent\textbf{Affine domains.}
For every $\op{x} = \SS(\textsf{x})$ in $\ADS^{1,n}$,
we define the \emph{affine domain} \index{Affine domain} (also denoted by $U(\textsf{x})$):
$$U(\op{x}) := \{ \op{y} \in \ADS^{1,n} \;/\; \langle \op{x} \mid \op{y} \rangle < 0 \}.$$

In other words, $U(\op{x})$ is the connected component of
$\ADS^{1,n}\setminus\SS(\textsf{x}^{\orth})$ containing $\op{x}$.

The boundary $\SS(\textsf{x}^{\orth}) \cap \ADS^{1,n}$ of $U(\op{x})$ in $\ADS^{1,n}$
has two components that are totally geodesic copies of the (Klein model of) hyperbolic space $\HH^n$.
One distinguish the \textit{past component} $H^-(\op{x})$
and the \textit{future component} $H^+(\op{x})$ characterized
by the following property: future oriented timelike geodesics enter $U(\op{x})$ through $H^-(\op{x})$
and exit through $H^+(\op{x})$. They are also called \textit{hyperplanes dual to $\op{x}$,} and we distinguish the hyperplane past-dual $H^-(\op{x})$ from the hyperplane future-dual $H^+(\op{x})$.

Every affine domain $U(\op{x})$, seen as a subset of $\Ein^{1,n}$, and lifted in
$\wt\AdS^{1,n} \approx \DD^n \times \RR$, is isometric to a region $\DD^{n} \times ]-\pi/2, +\pi/2[$.
\medskip

\noindent\textbf{Achronal subsets.} The description of achronal subsets of $\wt\AdS^{1,n}$ follows easily
from the description of achronal subsets of $\wt\Ein^{1,n}$: they are the graphs of $1$-Lipschitz maps
from $\DD^n$ into $\RR$. They are all contained in some affine domain as defined just previously.
\medskip

\noindent\textbf{Globally hyperbolic domains.}
Let $S$ be a closed edgeless achronal subset of $\wt\AdS^{1,n}$. In the conformal model, $S$ is the graph of a $1$-Lipschitz map  $f_S: \DD^n \to \RR$: this map uniquely extends to the boundary $\partial\DD^n$. The graph of this extension is a closed achronal edgeless subset $\partial S$ of $\partial\wt\AdS^{1,n}$, and the
Cauchy development dev$(S)$ is the intersection $\wt\AdS^{1,n} \cap E(\partial S)$ where $E(\partial S)$
is the invisible domain of $\partial S$ in $\wt\Ein^{1,n}$. It follows that dev$(S)$ is contained
in an affine domain too. In particular, dev$(S)$ is \emph{never} the entire
 anti-de Sitter space: this is another proof that $\wt\AdS^{1,n}$ is \emph{not} globally hyperbolic.
 Observe that $\partial S$ is purely lightlike if and only if $f_S(\bar{x}) = d_n(\bar{x,} \bar{x}_\infty)+t_0$
 where $t_0$ is a real number and $\bar{x}_\infty$ some point in $\partial\DD^n \subset \SS^n$.
\medskip

\noindent\textbf{The space of timelike geodesics.}
The content of this paragraph is mostly extracted from \cite[Section $4.5$]{quasifuchs}.
Timelike geodesics in $\AdS^{1,n}$ are intersections between $\AdS^{1,n} \subset \RR^{2,n}$ and $2$-planes $P$ in $\RR^{2,n}$
such that the restriction of $\mathrm{q}_{2,n}$ to $P$ is negative definite. The action of $\SO_0(2,n)$ on negative $2$-planes is transitive, and the stabilizer of the $(u, v)$-plane is $\SO(2) \times \SO(n)$. Therefore, the space of timelike geodesics is the symmetric space:
$$\cT_{2n} := \SO_0(2,n)/\SO(2) \times \SO(n).$$

\begin{rema}\label{rk:geodinvisible}
  Let $S$ be a  closed edgeless achronal subset of $\AdS^{1,n}$ such that $\partial S \subset \Ein^{1,n}$ is not purely lightlike.
Then, every timelike geodesic of $\AdS^{1,n}$ intersects $E(\partial S)$ (cf. Lemma 3.5 in \cite{merigot}), and since $E(\partial S)$ is convex, this intersection is connected, \ie is a single inextendible timelike geodesic of $E(\partial S)$. In other words,
one can consider $\cT_{2n}$ as the space of timelike geodesics of $E(\partial S)$ for any CEA $S$.
\end{rema}

$\cT_{2n}$ has dimension $2n$. We equip it with the Riemannian metric $g_\cT$ induced by the Killing form of $\SO_0(2,n)$. It is well known
that $\cT_{2n}$ has nonpositive curvature, and rank $2$: the maximal flats (\ie totally geodesic embedded Euclidean subspaces) have dimension $2$.
It is also naturally Hermitian. More precisely: let $\mathcal G = \mathfrak{so}(2,n)$ be the Lie  algebra of $G=\SO_0(2,n)$, and let $\cK$ be the
Lie algebra of the maximal compact subgroup $K:=\SO(2) \times \SO(n)$. We have the Cartan decomposition:
$$\mathcal G = \cK \oplus \cK^\perp$$
where $\cK^\perp$ is the orthogonal of $\cK$ for the Killing form. Then, $\cK^\perp$ is naturally identified with the tangent space at the origin of $G/K$. The adjoint action of the $\SO(2)$ term in the stabilizer defines a $K$-invariant complex structure on $\cK^\perp \approx T_K(G/K)$ that propagates through left translations to a genuine complex structure $J$ on $\cT_{2n} = G/K$. Therefore, $\cT_{2n}$
is naturally equipped with a structure of $n$-dimensional complex manifold, together with a $J$-invariant Riemannian metric, \ie a K\"ahler structure.

Here our purpose is to give another way to define this K\"ahler structure, starting from the anti-de Sitter space and clarifying the geometric nature of the associated symplectic form.

Let $\mathcal U^{1,n}$ be the space of future oriented timelike tangent vectors of norm $-1$ of $\AdS^{1,n}$.
Since the tangent space of $\AdS^{1,n}$ at a point $\textsf{x}$ is naturally identified with $\textsf{x}^\perp$,
there is a natural identification between $\mathcal U^{1,n}$ and the pairs $(\textsf{x,} \textsf{y})$ of elements
of $\AdS^{1,n}$ satisfying $\textsf{y} \in H^+(\textsf{x})$. The tangent bundle of $\mathcal U^{1,n}$ at a point
$(\textsf{x,} \textsf{y})$ is naturally identified with elements $(\dot{\textsf{x},} \dot{\textsf{y}})$ of $\RR^{2,n} \times \RR^{2,n}$ satisfying:
\begin{eqnarray*}
  \sCa{\dot{\textsf{x}}}{\textsf{x}} &=& 0, \\
  \sCa{\dot{\textsf{y}}}{\textsf{y}} &=& 0, \\
  \sCa{\dot{\textsf{x}}}{\textsf{y}} + \sCa{\dot{\textsf{y}}}{\textsf{x}} &=& 0.
\end{eqnarray*}

There is a canonical $\SO_0(2,n)$-invariant pseudo-Riemannian metric $\| . \|$ on $\mathcal U^{1,n}$:
$$\| (\dot{\textsf{x},} \dot{\textsf{y}}) \| := q_{2,n}(\dot{\textsf{x}}) + q_{2,n}(\dot{\textsf{y}}).$$

A quick computation shows that $\| . \|$ is Lorentzian. Moreover, it is preserved by the
geodesic flow $\Phi^t$ on $\mathcal U^{1,n}$ that can be defined by:
$$\Phi^t(\textsf{x,} \textsf{y}) = (\cos(t)\textsf{x} + \sin(t)\textsf{y,} -\sin(t)\textsf{x} + \cos(t)\textsf{y}).$$
The Killing vector field generating $\Phi^t$ is $Z(\textsf{x,} \textsf{y}) = (\textsf{y,} -\textsf{x})$, hence
of $\| . \|$-norm $-2$: the $\| . \|$-orthogonal $Z^\perp$ is therefore spacelike.

Now $\cT_{2n}$ is naturally identified with the orbit space of $\Phi^t$: we write elements
of $\cT_{2n}$ as equivalence classes $[\textsf{x,} \textsf{y}]$ of the orbital relation on $\mathcal U^{1,n}$.
We have a natural way to define an $\SO_0(2,n)$-invariant Riemannian metric $\overline{\|} . \overline{\|}$ on $\cT_{2n}$:
for every $[\textsf{x,} \textsf{y}]$ and every tangent vector $\xi$ at $[\textsf{x,} \textsf{y}]$, select
a representative $(\textsf{x,} \textsf{y})$ of $[\textsf{x,} \textsf{y}]$ and the unique vector $v$ tangent
to $\cT_{2n}$, orthogonal to $Z(\textsf{x,} \textsf{y})$, and projecting onto $\xi$. Define $\overline{\|}(\xi)\overline{\|}$ as
$\|v\|$: since $\Phi$ is isometric, this value does not depend on the choice of $(\textsf{x,} \textsf{y})$. This
defines a Riemannian metric on $\cT_{2n}$ that happens to be (up to a constant factor) the metric on
the symmetric space $\cT_{2n} := \SO_0(2,n)/\SO(2) \times \SO(n)$.

The symplectic form $\omega$ associated with the K\"ahler-Hermitian structure of $\cT_{2n}$ can be described as follows: the Liouville $1$-form $\lambda$ on $\mathcal U^{1,n}$ is defined by:
$$\lambda_{(\textsf{x,} \textsf{y})}(\dot{\textsf{x},} \dot{\textsf{y}}) = -\sCa{\dot{\textsf{x}}}{\textsf{y}} = \sCa{\dot{\textsf{y}}}{\textsf{x}}.$$
The contact hyperplane, kernel of $\lambda$, is the orthogonal $Z^\perp$.
The geodesic flow $\Phi^t$ is actually the Reeb flow for $\lambda$: we have $\lambda(Z) = +1$ and $L_Z\lambda = 0,$ so
the exterior derivative $\hat{\omega} = d\lambda$ is $\Phi^t$-invariant and for any $v$ in $Z^\perp$ we have
$\hat{\omega}(Z, v) = i_Zd\lambda(v) = i_Zd\lambda(v) + di_Z\lambda = L_Z\lambda = 0$. Hence we can define on
$\cT_{2n}$ the $2$-form $\omega(\xi, \xi') = \hat{\omega}(v,v')$ where $v$, $v'$ are the lifts of $\xi$, $\xi'$ above $v$, $v'$ orthogonal to $Z$. This $2$-form is closed and non-degenerate, \ie a symplectic form on $\cT_{2n}$.

Finally, we define the complex structure: first observe that there is a natural complex structure on
each contact hyperplane $Z(\textsf{x,} \textsf{y})^\perp$: this hyperplane is identified with the set of pairs
$(\dot{\textsf{x},} \dot{\textsf{y}})$ with $\dot{\textsf{x}}$ and $\dot{\textsf{y}}$ both in $\textsf{x}^\perp \cap \textsf{y}^\perp$. Then $(\dot{\textsf{x},} \dot{\textsf{y}}) \mapsto (\dot{\textsf{y},} -\dot{\textsf{x}})$
is an involution in $Z(\textsf{x,} \textsf{y})^\perp$, commuting with the geodesic flow, and induces the complex structure $J$ on $\cT_{2n}$.
Observe that $J$ satisfies $\omega(v, J(v)) = \| v \|$: it is a calibration between $\omega$ and $\| . \|$.

\medskip

\noindent\textbf{Gauss map.}
Let $S$ be a differentiable Cauchy hypersurface in a GH domain $E(\partial S)$. The \textit{Gauss map} \index{Gauss map} of $S$
is the map $\nu: S \to \cT_{2n}(\rho)$ that sends every element $x$ of $S$ to the unique timelike geodesic orthogonal to $S$ at $x$.

Since every timelike geodesic intersects $S$ at most once, the Gauss map is always injective. The image of the Gauss map
is actually the set of timelike geodesics that are orthogonal to $S$. Since every timelike geodesic intersects $S$, it follows easily that
the image of the Gauss map is closed, and that the Gauss map is actually an embedding.

Assume that $S$ is the image of a smooth spacelike embedding $\textsf{x}: \Sigma \to \AdS^{1,n}$. Then
it induces a map $\textsf{y}: \Sigma \to \AdS^{1,n}$ where $\textsf{y}(p)$ is the dual of the unique totally
geodesic hypersurface tangent to $S$ at $\textsf{x}(p)$, and $\nu(p) = [\textsf{x}(p), \textsf{y}(p)]$.
Here we point out that the linear map $d{\textsf{x}}(\dot{p}) \to d{\textsf{y}}(\dot{p})$ is the shape operator $B$ at $\textsf{x}(p)$,
so that the second fundamental form is II$(\dot{p}) = \langle d{\textsf{x}}(\dot{p}) \mid B(d{\textsf{x}}(\dot{p})) \rangle = \langle d{\textsf{x}}(\dot{p}) \mid d{\textsf{y}}(\dot{p}) \rangle$.


We also point out the following consequence of the equality $\langle d{\textsf{x}}(\dot{p}) \mid {\textsf{y}}({p})) \rangle =0:$ \emph{the image of the Gauss map is a Lagrangian submanifold of $\omega$.} Conversely,
let $\varphi: \Sigma \to \cT_{2n}$ be a Lagrangian immersion for some $n$-dimensional simply connected manifold $\Sigma$. One can lift $\varphi$ to some immersion $\hat{\varphi}: \Sigma \to \mathcal U^{1,n}$ orthogonal to $Z$: indeed, select a base point $p_0$ in $\Sigma$, and for any $p$ in $\Sigma$ let $\alpha: [0, 1] \to \Sigma$ and $\beta: [0, 1] \to \Sigma$
be two path with $\alpha(0)=\beta(0)=p_0$ and $\alpha(1)=\beta(1)=p$. Then, lift them to paths $\hat{\alpha}$, $\hat{\beta}$ in $\mathcal U^{1,n}$ orthogonal to $Z$ and such that $[\hat{\alpha}(t)] = \varphi(\alpha(t))$ and
$[\hat{\beta}(t)] = \varphi(\beta(t))$. Assume $\hat{\alpha}(1)$ and $\hat{\beta}(1)$ are both above $\alpha(1)=\beta(1)=p$, hence there is
a real number $t$ such that $\hat{\beta}(1) = \Phi^t(\hat{\alpha}(1))$. The loop obtained by composing $\hat{\alpha}$, the portion of the $\Phi$-orbit between $\hat{\alpha}(1)$ and $\hat{\beta}(1)$ and the inverse of $\hat{\beta}$ is homotopically trivial, hence the boundary
of a disk $D$. Then, since the integral of $\lambda$ along the portion of $\Phi$-orbit is $t$, we have:
$$\int_D \varphi^\ast\omega = \int_{\hat{\alpha}}\lambda +t - \int_{\hat{\beta}}\lambda.$$
Now observe that $\int_{\hat{\alpha}}\lambda = \int_{\hat{\beta}}\lambda = 0$, hence, since $\int_D \varphi^\ast\omega =0$ (because $\varphi$ is Lagrangian) we get $t=0$. The equality $\hat{\alpha}(1)=\hat{\beta}(1)$ follows: we can define $\hat{\varphi}(p) = \hat{\alpha}(1)$.

The immersion $\hat{\varphi}: \Sigma \to \mathcal U^{1,n}$ is spacelike for $\|$ (since it is orthogonal to $Z$). Write $\hat{\varphi}(p) = (\textsf{x}(p), \textsf{y}(p))$. It may happen that $p \mapsto \textsf{x}(p)$ is not an immersion at some point $p_0$, but then, near $p_0$, simply replace $\hat{\varphi}$ by $\hat{\varphi}_t = \Phi^t \circ \hat{\varphi}$ by any non zero real number $t$ (it amounts to replacing $\textsf{x}(p)$ by $\cos(t)\textsf{x}(p) + \sin(t)\textsf{y}(p)$): after this modification, $p \mapsto \textsf{x}(p)$ becomes a spacelike immersion near $p_0$, whose Gauss map is the restriction of $\varphi$ near $p_0$. Of course, this construction is local and we may fail to find a $t$ valid over the entire $\Sigma$.

\section{The classification of MGHC spacetimes}
\label{sec:classGH}
In this section, we present the classification of MGHC spacetimes with constant curvature. This classification
has been initiated by G. Mess who introduced many ideas (\cite{mess1}), and continued by several authors (\cite{barflat, bonsante, scannell}). This classification is almost complete, but there are some remaining questions.

\subsection{The flat case}\label{sub:flatGH}
The simplest example of MGHC flat spacetime is the quotient of the entire Min\-kowski space by a discrete group $\Gamma$ of translations, such that the translation vectors of $\Gamma$ form a lattice in $\RR^{0,n}$. These
spacetimes are geodesically complete. They are isometric to the product $(\mathbb T^n, \bar{h}_0) \times (\RR, -dt^2)$, where $(\mathbb T^n, \bar{h}_0)$ is a flat torus. We call them \emph{translation spacetimes.}\index{translation spacetime}

The most fundamental example of MGHC flat spacetime is the quotient of the cone $I^+_0$ in $\RR^{1,n}$ by a cocompact lattice $\Gamma$ of $\SO_0(1,n)$:
$$M_0(\Gamma) = \Gamma\backslash I^+_0.$$
The restriction of $-q_{1,n}$ to $I^+_0$ is a time function whose levels sets are Cauchy hypersurfaces --- in particular,
the hyperboloid model $\HH^n$ of the hyperbolic space which is the $1$-level set.

We call this example the \emph{standard conformally static example}.

Identify $\RR^{1,n}$ with the domain Mink$^+(x_0)$ of $\wt\Ein^{1,n}$.
The closure of $\HH^n$ in $\wt\Ein^{1,n}$ is then the union of $\HH^n$ with a closed edgeless achronal set $\partial\HH^n \subseteq \mathcal I^+$. Then, $I^+_0$ coincides with the invisible domain
$E(\partial\HH^n) = \wt\Ein^{1,n} \setminus \left(J^-(\partial\HH^n)\cup J^+(\partial\HH^n)\right)$.

More generally, let $\Lambda$ be a (non edgeless) closed achronal subset of $\wt\Ein^{1,n}$ contained in the conformal boundary
$\mathcal I^+$ of a Minkowski domain Mink$^+(x_0) \approx \RR^{1,n}$, invariant by a
 torsionfree discrete subgroup $\Gamma$ of Isom$(\RR^{1,n})$.
Then, $\Omega(\Lambda) = E(\Lambda) \cap$ Mink$^+(x_0)$ is globally hyperbolic, $\Gamma$-invariant, and the action of $\Gamma$ on $\Omega(\Lambda)$
is free and proper. The quotient $M_\Lambda(\Gamma)$ is a globally hyperbolic spacetime, but not necessarily Cauchy compact. The cosmological time of $\Omega(\Lambda)$ --- that is the lift of the cosmological time of $M_\Lambda(\Gamma)$ --- is regular. Moreover,
if $\Lambda$ contains no proper $\Gamma$-invariant closed subset, then $M_\Lambda(\Gamma)$ is maximal among
flat GH spacetimes. Finally, $M_\Lambda(\Gamma)$ is future complete, in the sense that any future oriented timelike ray is geodesically complete.

We have observed that elements of $\mathcal I^+$ correspond to past-half spaces in $\RR^{1,n}$, and $\Omega(\Lambda)$ is obtained by removing all the half-spaces corresponding to elements of $\Lambda$. It follows that $\Omega(\Lambda)$ is an
intersection of half-spaces, hence convex. We recover the notion of \emph{regular domain} as defined in \cite{bonsante}. The boundary of $\Omega(\Lambda)$ in $\RR^{1,n}$ is very interesting: it is a CEA, that is, the graph of a $1$-Lipschitz map $f$ whose differential has norm $1$
almost everywhere.

Observe that if we replace $\mathcal I^+$ by $\mathcal I^-$, the result will be a similar spacetime, but that is past complete and not future complete.

The following Theorem was proved by Mess in dimension $2+1$ (\cite{mess1}), by Bonsante in the case where the holonomy is assumed to have a discrete cocompact linear part in $\SO_0(1,n)$ (\cite{bonsante}), and independently and in full generality in \cite{barflat}:

\begin{teo}\
  Up to finite coverings, every MGHC flat spacetime is isometric to either a translation spacetime, or a quotient $M_\Lambda(\Gamma)$, where $\Gamma$ is a discrete subgroup of \/$\op{Isom}(\RR^{1,2})$.
\end{teo}

We want to describe further $M_\Lambda(\Gamma)$ when it is Cauchy compact. Most of the following claims are non-trivial and we refer to \cite{barflat} for their proofs. Notice first that in this case, the achronal subset $\Lambda$ of $\mathcal I^\pm \approx \RR \times \SS^{n-1}$ is necessarily a topological sphere, the graph
of a continuous map $f: \SS^{n-1} \to \RR$. The invisible domain $E(\Lambda)$ is then contained in Min$_+(x_0)$, hence we have $\Omega(\Lambda) = E(\Lambda)$.

Consider first the \textit{proper} case, \ie the case where the closure of $\Omega(\Lambda)$ does not contain any affine line of
$\RR^{1,n}$:  $\Gamma$ is then isomorphic to a cocompact lattice of $\SO_0(1,n)$. More precisely, the linear part
homomorphism Isom$(\RR^{1,n}) \to \SO_0(1,n)$ is faithful and has discrete and cocompact image. In this case, we call $M_\Lambda(\Gamma)$ a \emph{standard spacetime.} \index{standart spacetime}

In general, excluding the particular case of \emph{Misner spacetimes,} \index{Misner spacetime} up to finite coverings, a MGHC spacetime is a twisted product over a standard spacetime by flat tori. In particular, Cauchy hypersurfaces are always finite covers of products of hyperbolic closed manifolds by tori.
See \cite{barflat} for more details.

There is also an interesting case of flat MGH spacetime, but not Cauchy compact: the \emph{unipotent spacetimes.} \index{unipotent spacetime}
A unipotent spacetime is the quotient of either a half-space in $\RR^{1,n}$ bounded by a lightlike hyperplane, or the region between two parallel lightlike hyperplanes, by a discrete nilpotent group $\Gamma$ whose linear part is a discrete subgroup of the stabilizer in $\SO_0(1,n)$
of a point of $\partial\HH^n$. For more details, see \cite{barflat}.

\subsection{The \mathversion{bold}$\dS$ case}\label{sub:dSGH}
In some way, MGHC Cauchy compact spacetimes locally modeled on $\dS^{1,n}$ first appeared in a paper
by Kulkarni and Pinkall (\S $3.4$ of \cite{kulkarni}), but the authors did not insist on the de Sitter nature of
the spaces they were considering, and, presumably, were not aware of their interpretation as globally hyperbolic spacetimes.
The fact that these examples give the complete list of MGHC de Sitter spacetimes was proved by K. Scannell (\cite{scannell}), involving some ideas of Mess. There is also an alternative description in \cite{ABBZ} which is the one we use here.

Unlike the other two cases (the flat and anti-de Sitter cases), locally de Sitter MGHC spacetimes are not in general quotients of
open domains in $\dS^{1,n}$.

The crucial point is that locally de Sitter MGHC spacetimes are in one-to-one correspondence with \emph{closed M\"obius manifolds,} \index{M\"obius manifold}  \ie closed manifolds locally modeled on the conformal sphere $\SS^n$. This correspondence involves the identification of $\dS^{1,n}$ with the space of round discs in $\SS^n$. It goes as follows: let $\Sigma$ be a closed manifold of dimension $n$, locally modeled on $(\SS^n, \SO_0(1,n+1))$. Let $\cD: \wt\Sigma \to \SS^n$ be the developing map of this
$(\SS^n, \SO_0(1,n+1))$-structure, and let $\rho: \Gamma = \pi_1(\Sigma) \to \SO_0(1,n+1)$ be the holonomy representation. Let $\mathcal B(\wt\Sigma)$ be the space of open domains $B$ in $\wt\Sigma$
such that the restriction of $\cD$ to $B$ is a homeomorphism onto a round disc of $\SS^n$. Then the action of $\Gamma$ on $\mathcal B(\wt\Sigma)$ is free and proper: let $\mathcal B(\Sigma)$ be the quotient space. The developing map $\cD$ induces a $\Gamma$-equivariant local homeomorphism $\widehat\cD: \mathcal B(\wt\Sigma) \to \mathcal B(\SS^n) \approx \dS^{1,n}$, hence a locally de Sitter structure on $\mathcal B(\Sigma)$. As a locally de Sitter
manifold, $\mathcal B(\Sigma)$ is maximal globally hyperbolic.

Actually, we observed that there are two ways to identify $\dS^{1,n}$ with $\mathcal B(\SS^n)$: for one of them the spacetime $\mathcal B(\Sigma)$ is geodesically complete in the future, and for the other one, $\mathcal B(\Sigma)$ is geodesically complete in the past.

\begin{teo}[K. Scannell \cite{scannell}]
Every maximal globally hyperbolic Cauchy compact locally de Sitter spacetime is isometric to
the spacetime $M(\Sigma)$ associated with a $(\SS^n, \SO_0(1,n+1))$-manifold $\Sigma$.
\end{teo}

Thurston observed  that $(\SS^2, \SO_0(1,3))$-manifolds are in one-to-one correspondence with \emph{hyperbolic ends,} \ie hyperbolic $3$-manifolds homeomorphic to $\Sigma \times [0, +\infty[$ with a concave boundary $\Sigma \times \{ \; 0 \}$ and complete at the end $\Sigma \times \{ +\infty \}$ (for more details, see for example \cite{benbon}). This generalizes in any dimension: every locally de Sitter MGHC spacetime $M(\Sigma)$ has an associated dual hyperbolic manifold homeomorphic to $\Sigma \times [0, +\infty[$. See \cite{ABBZ} for more details.

We point out particular elementary cases:

-- \emph{the elliptic case}: this is the case where $\Sigma$ is the round sphere $\SS^n$; $M(\Sigma)$ is then
the de Sitter space;

-- \emph{the parabolic case}: this is the case where $\Sigma$ is a quotient of the flat conformal Euclidean space $\RR^n$, \ie the once punctured sphere. Then, in the geodesically future complete case, $M(\Sigma)$ is a quotient of the complement in $\dS^{1,n}$ of the past of a point in the future conformal boundary. It is dual to the hyperbolic end corresponding to one hyperbolic cusp. In particular, Cauchy surfaces are finite quotients of umbilic tori.

In the remaining non-elementary \emph{hyperbolic case,} $M(\Gamma)$ has a regular cosmological time.

\subsection{The \mathversion{bold}$\AdS$ case}
\label{sec.regads}
A convenient recent reference for the content of this section is \cite{quasifuchs}.
Let $\wt\Lambda$ be a closed edgeless achronal subset of $\partial\wt\AdS^{1,n} \approx \wt\Ein^{1,n-1}$ ($n \geq 2)$.
Consider it as a closed achronal subset of $\wt\Ein^{1,n}$. Assume that it is not purely lightlike. Then, the invisible domain $E(\wt\Lambda)$ in $\wt\Ein^{1,n}$ has two connected components: one contained in $\wt\AdS^{1,n}$ and the other in the second anti-de Sitter component of $\wt\Ein^{1,n} \setminus \wt\Ein^{1,n-1}$. We call $\Omega(\wt\Lambda)$ the first component.
We can write:
$$\Omega(\wt\Lambda)=:\wt\AdS^{1,n} \setminus
\left(J^-(\wt\Lambda)\cup J^+(\wt\Lambda)\right)
$$

We denote by 
$\Omega(\Lambda)$ the projection of $\Omega(\wt \Lambda)$ in $\AdS^{1,n}$ (cf. Lemma \ref{le.one-to-one}).

The domains $\Omega(\Lambda)$ and $\Omega(\wt\Lambda)$ are isometric to each other. They have regular cosmological time, in particular, they are GH.

Let $\Gamma$ be a torsionfree discrete subgroup of $\SO_0(2,n)$ preserving $\Lambda$. Then, the action of
$\Gamma$ on $\Omega(\Lambda)$ is free and properly discontinuous, and preserves the cosmological time $\tau$.
The quotient $M_\Lambda(\Gamma) := \Gamma\backslash\Omega(\Lambda)$ is a MGH spacetime locally modeled on the anti-de Sitter space. We call it a \emph{regular MGH anti-de Sitter spacetime.} \index{regular MGH anti-de Sitter spacetime}

\begin{teo}
Every locally anti-de Sitter MGHC spacetime is isometric to a regular MGHC AdS spacetime.
\end{teo}


We now give a more detailed presentation of the geometric features of these spacetimes.
In the conformal model, $\Omega(\wt\Lambda)$ is the region in $\DD^n \times \RR$ between the graphs
of two $1$-Lipschitz maps $f^\pm: \DD^n \to \RR$ that are extensions of the map $f_\Lambda: \partial\DD^n \to \RR$ whose graph is $\Lambda$.
The graph of $f^-$ (respectively $f^+$) is a closed achro\-nal subset of $\wt\AdS^{1,n}$ that we call the \textit{lifted past}
(respectively \textit{future}) \textit{horizon} of $\Omega(\wt\Lambda)$, and denote by $\wt\cH^-(\Lambda)$ (respectively $\wt\cH^+(\Lambda)$).
The projections in $\AdS^{1,n}$ of $\wt\cH^\pm(\Lambda)$ are called \textit{past} and \textit{future horizons} of $\Omega(\Lambda)$, and denoted by $\cH^\pm(\Lambda)$.

Consider now $\Lambda$ as a closed subset in $\SS(\mathcal C_{n}) \subset \SS(\RR^{2,n})$, boundary of the Klein model $\ADS^{1,n}$: since $\Lambda$ is achronal, the convex hull $\op{Conv}(\Lambda)$ is contained in the closure of $\ADS^{1,n}$. Actually, its intersection with $\SS(\mathcal C_{n})$ is precisely
$\Lambda$; in parti\-cular,
$E(\Lambda)$ characterizes $\Lambda$.
It happens (\cite{ABBZ}) that
$E(\Lambda)$ is the interior of the dual $\op{Conv}(\Lambda)^\ast$ --- in particular, $E(\Lambda)$
contains the interior of $\op{Conv}(\Lambda)^\ast$.
\medskip

\noindent\textbf{The Fuchsian case.}
There is a particular case: the case where $\Lambda$ is the boundary of a totally geodesic copy of $\HH^n$ in $\AdS^{1,n}$. THis means
that $\Lambda$ is the graph of an affine map. Then, $\op{Conv}(\Lambda)$ has empty interior, and the complement of $\Lambda$ in $\op{Conv}(\Lambda)$
is the totally geodesic subspace bounded by $\Lambda$. The group $\Lambda$ is then a discrete subgroup of a conjugate of $\SO_0(1,n)$ in $\SO_0(2,n)$. We call this case the \textit{Fuchsian case.}\index{Fuchsian AdS spacetime}
\medskip

\noindent\textbf{Past tight region.}
From now on, we assume that $\Lambda$ is not Fuchsian. Its complement in the boundary $\partial\op{Conv}(\Lambda)$ has
two connected components. Both are closed achronal subsets of $\AdS^{1,n}$. More precisely,
in the conformal model their lifs to $\wt{\AdS}^{1,n}$ are graphs
of 1-Lipschitz maps $F^+$, $F^-$ from $\DD^n$
into $\RR$ such that
\begin{equation}
\label{eq:fF}
f^- \leq F^- \leq F^+ \leq f^+.
\end{equation}

The graph of $F^-$ is the \textit{past component}
$S^-(\Lambda)$ and the graph of $F^-$ is the \textit{future component} $S^+(\Lambda)$. The region between the past horizon $\cH^-(\Lambda)$ and the future component $S^+(\Lambda)$ is the
\emph{past tight region} \index{past tight region} and it is denoted by $E_0^-(\Lambda)$.

Since $E(\Lambda)$ and $\op{Conv}(\Lambda)$ are convex and dual to each other, for every element
$x$ in $S^-(\Lambda)$
(respectively $S^+(\Lambda)$) there is an element $p$ of $\Lambda$ or $\cH^+(\Lambda)$ (respectively $\cH^-(\Lambda)$) such that $H^-(p)$
(respectively $H^+(p)$) is a support hyperplane for $S^-(\Lambda)$
(respectively $S^+(\Lambda)$) at $x$: these support hyperplanes are either totally geodesic copies of $\HH^n$ (if $p \in \AdS^{1,n}$) or degenerate (if $p \in \Lambda$).

\begin{figure}[ht]
  \centering
  \includegraphics[width=8cm]{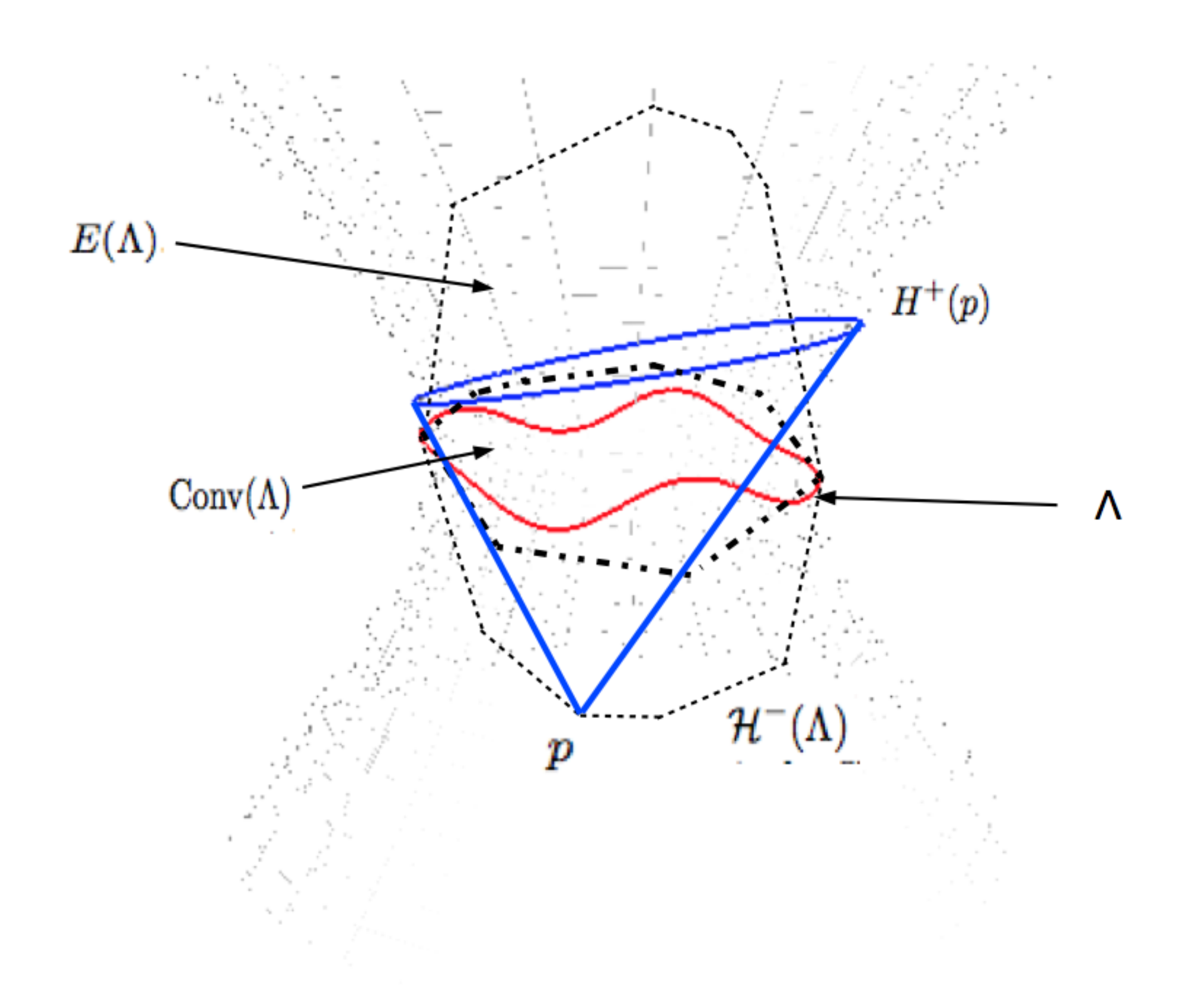}
\centering  \caption{The global situation. The dotted hyperboloid represents the boundary of an affine domain of
$\AdS^{1,n}$ containing the invisible domain $E(\Lambda)$. The limset $\Lambda$ is represented by a topological circle turning around
the hyperboloid, and $\op{Conv}(\Lambda)^\circ$ is a convex subset inside the (dual) convex subset $E(\Lambda)$. The
future-dual plane $H^+(p)$ for $p$ in the past horizon $\cH^-(\Lambda)$ is a support hyperplane of $S^+(\Lambda)$.}
  \label{fig:situa}
\end{figure}

Similarly, at every element $x$ of $\cH^-(\Lambda)$ (respectively $\cH^+(\Lambda)$) there is a support hyperplane $H^-(p)$ (respectively $H^+(p)$) where $p$ is an element of
$S^+(\Lambda) \cup \Lambda$ (respectively $S^-(\Lambda) \cup \Lambda$). See Figure \ref{fig:situa}.

\medskip

\noindent\textbf{Cosmological lines.} The past tight region $E_0^-(\Lambda)$ is precisely the region
where the cosmological time $\tau$ takes value $<\pi/2$:
$$E_0^-(\Lambda) = \{ \tau < \pi/2 \}.$$
For every $x$ in $E_0^-(\Lambda)$ there is a unique \emph{realizing geodesic} for $x$. More
precisely, there is one and only one element $r(x)$ in the past horizon $\cH^-(\Lambda)$ --- called the \emph{cosmological retract of $x$} \index{cosmological retract} ---
such that the segment $]r(x), x]$ is a timelike geodesic whose Lorentzian length is precisely the cosmological time $\tau(x)$.

The restriction of $\tau$ to $E_0^-(\Lambda)$ is $C^{1,1}$, \ie it is differentiable with locally Lipschitz derivative) (see \cite{mess1, benbon} for the case $n=2$, and \cite{quasifuchs} for the general case).
The realizing geodesics are orthogonal to the level sets of $\tau$. The space of realizing
geodesics is then an $n$-dimensional closed embedded Lipschitz submanifold.
We denote it by $\Sigma(\Lambda)$ and call it the \textit{space of cosmological lines.} \index{space of cosmological lines}
\medskip

\noindent\textbf{The initial singularity.}
For every $p$ in $\cH^-(\Lambda)$, $H^+(p)$ is a support hyperplane for $\op{Conv}(\Lambda)$, but it could be at a point in $\Lambda$.
Elements of $\cH^-(\Lambda)$ that are dual to support hyperplanes for $\op{Conv}(\Lambda)$ at a point inside $\AdS^{1,n}$, \ie in $S^+(\Lambda)$,
form an interesting subset of $\cH^-(\Lambda)$, the \textit{initial singularity set} \index{initial singularity} (cf. \cite{benbon}).
\medskip

\noindent\textbf{Split \mathversion{bold}$\AdS$-spacetimes.}
Consider the sum $\RR^{1,k} \oplus \RR^{1,\ell}$ with $k+\ell = n$, equipped with the quadratic form
$q_{1,k} + q_{1,\ell}$. It is isometric to $(\RR^{2,n,} q_{2,n})$. It provides an embedding
of $\SO_0(1,k) \times \SO_0(1,\ell)$ in $\SO_0(2,n)$. It preserves in $\partial\AdS^{1,n}$ an achronal topological sphere $\Lambda_{k, \ell}$, link of two spheres $\SS^{k-1}$ and $\SS^{\ell-1}$, where every point
in $\SS^{k-1}$ is linked to every point in $\SS^{\ell-1}$ by the unique future oriented lightlike segment in $\partial\AdS^{1,n}$ going from the point in $\SS^{k-1}$ towards the point in $\SS^{\ell-1}$.
The associated domain $E(\Lambda_{k, \ell})$ is then globally hyperbolic, and for every cocompact lattice $\Gamma$ in $\SO_0(1,k) \times \SO_0(1,\ell)$ the quotient $\Gamma\backslash E(\Lambda_{k, \ell})$ is a MGHC $\AdS$-spacetime, called a \emph{ split $\AdS$-spacetime.} \index{split AdS spacetime}
Observe that the achronal sphere $\Lambda_{k, \ell}$ is not acausal. A natural conjecture is that
any MGHC $\AdS$-spacetimes either is a split $\AdS$-spacetime, or has an associated achronal subset $\Lambda$ that
is acausal. For more details, see \cite{franceseinstein} or \cite{quasifuchs}.

\subsection{Invisible domains}
\label{sub:btz}
Let $\Lambda$ be a non-purely lightlike closed achronal subset of $\partial\AdS^{1,n} \approx \Ein^{1,n-1}$, but not necessarily edgeless: $\Lambda$ is the graph of a $1$-Lipschitz map $f: \Lambda_0 \to \RR$ where $\Lambda_0$ is a closed subset of $\SS^{n-1}$. Assume that it is preserved by a torsionfree discrete subgroup $\Gamma$ of $\SO_0(2,n)$. We consider it as a subset of $\Ein^{1,n}$: let $E(\Lambda)$ be the invisible domain in $\Ein^{1,n}$, $D(\Lambda)$ be the restriction of $E(\Lambda)$ to $\Ein^{1,n-1}$ and $\Omega(\Lambda)$ be its restriction to $\AdS^{1,n}$. Then, the quotient $M(\Lambda) := \Gamma\backslash\Omega(\Lambda)$ is a strongly causal spacetime, diffeomorphic to a product $S \times \RR$, but not globally hyperbolic if $\Lambda_0 \neq \SS^{n-1}$: some causal curves may escape from it by some point in $\Gamma\backslash D(\Lambda)$.

Actually, one can consider $\Gamma\backslash D(\Lambda)$ as the conformal boundary of $M(\Lambda)$.
The past of $\Gamma\backslash D(\Lambda)$ in $M(\Lambda)$ is the quotient of $J^-(D(\Lambda)) \cap \Omega(\Lambda)$ by $\Gamma$. Its complement in $M(\Lambda)$ is globally hyperbolic: it is the quotient of
$\Omega(\Lambda^+)$ by $\Gamma$, where $\Lambda^+$ is the future boundary of $D(\Lambda)$, \ie the graph of the maximal $1$-Lipschitz extension of $f$ on $\SS^{n-1}$.

In Section \ref{sub:btz2} we will give a geometrico-relativistic interpretation of this feature as
idealistic models of spacetimes containing black-holes.

\section{Discrete groups of isometries}
\label{sec:discrete}
In this section, we explore how the notion of global hyperbolicity may help to
understand actions of groups of isometries on the model spacetimes, in a way similar
to the traditional theory of Kleinian groups. We focus on the notion of \textit{achronal} subgroups,
and show how they allow to develop a theory similar to the classical theory of groups of
isometries of the hyperbolic space $\HH^n$, in particular, that such a subgroup has always a canonical
\emph{limit set} in the conformal boundary.

We then clarify the relation between achronal subgroups and \emph{Anosov representations} in the context of isometry groups of spacetimes of constant curvature.

\begin{defi}
Let $\Gamma$ be a group of isometries of a spacetime $(M, g)$. An orbit $\Gamma.p$ is \emph{achronal} \index{achronal orbit} if for every non-trivial element $\gamma$ of $\Gamma$ the iterate $\gamma p$ is not in the strict future or past $I^\pm(p)$ of $p$. The group $\Gamma$ itself is \emph{achronal} \index{achronal group} if it admits an achronal orbit.
\end{defi}

In some way, one may extend this notion to any pseudo-Riemannian metric, and it happens that
isometry groups of Riemannian spaces are always achronal since there is no timelike tangent vectors in this case.


\begin{remark}
We do not assume $\Gamma$ to be a discrete subgroup, even if it is the main case we have in mind. Discreteness is not necessary for this basic theory of limit sets.
\end{remark}

\subsection{The de Sitter case}
In this case $\Gamma$ is a subgroup of $\SO_0(1,n+1)$, hence can be seen as a group of isometries of $\HH^{n+1}$.
The key observation is that, in our framework, $\Gamma$ is always achronal.

Indeed, the boundary
of $\dS^{1,n}$ in $\wt\Ein^{1,n}$ is the union of two spacelike spheres: the past and future conformal
boundaries $\partial_\pm\dS^{1,n}$ which are achronal.
Actually, every component of $\partial\dS^{1,n}$ is conformally equivalent to $\SS^n \approx \partial\HH^{n+1}$,
and the classical theory of limit sets for groups of hyperbolic isometries provides a limit
set $\Lambda_\Gamma$, usually defined as the set of points in $\partial\HH^{n+1}$ that are accumulation points of any orbit in $\HH^{n+1}$.

Here, we have two copies of the limit set in $\partial\dS^{1,n}$: one $\Lambda^+_\Gamma$ in $\partial_+\dS^{1,n}$, and the other one $\Lambda^-_\Gamma$ in $\partial_-\dS^{1,n}$. Each of them
is acausal in $\wt\Ein^{1,n}$, hence we can define as before the (non-empty) regions  $\Omega_\pm(\Lambda_\Gamma)$, intersections
between $\dS^{1,n}$ and the invisible domain $E(\Lambda^\pm_\Gamma)$ in $\wt\Ein^{1,n}$. Observe that $\Omega_\pm(\Lambda_\Gamma)$ is globally hyperbolic, since $\dS^{1,n}$ and $E(\Lambda^\pm_\Gamma)$ are GH domains in $\wt\Ein^{1,n}$. The
union $\Omega(\Lambda_\Gamma)$ of $\Omega_+(\Lambda_\Gamma)$ and $\Omega_-(\Lambda_\Gamma)$ can be characterized as the interior of the set of points in $\dS^{1,n}$ with achronal $\Gamma$-orbit, and
$\Lambda_\Gamma = \Lambda^+_\Gamma \cup \Lambda^-_\Gamma$
is the set of non-trivial accumulation points of \textbf{achronal} orbits $\Gamma.p$ in $\dS^{1,n}$
$$\Lambda_\Gamma := \overline{\Gamma.p} \setminus \Gamma p, \mbox{ for any $p$ with achronal $\Gamma$-orbit.}$$

$\Lambda_\Gamma$ is empty if and only if $\Gamma$ is relatively compact in $\SO_0(1,n)$, in which case
one can take $\Omega(\emptyset)$ as the entire de Sitter space --- actually, in this case $\Gamma$ is conjugate to a subgroup of the maximal subgroup $\SO(n)$ of $\op{SO}_0(1,n)$: it preserves a foliation
by umbilical spacelike spheres: every orbit is acausal, and the action on $\dS^{1,n}$ is proper.

When $\Gamma$ is discrete, $\Omega(\Lambda_\Gamma)$ can be seen as the natural domain of $\dS^{1,n}$ on which
$\Gamma$ acts properly discontinuously.

\begin{remark}
  Consider the Klein model in $\mathbb P(\RR^{1,n+1})$: the convex hull $\op{Conv}(\Lambda_\Gamma)$
  is contained in the Klein model $\mathbb  P(\{ q_{1,n+1} < 0\})$ of the hyperbolic space, and $\Omega(\Lambda_\Gamma)$ is the intersection
  between the Klein model $\mathbb  P(\{ q_{1,n+1} > 0\})$ of the de Sitter space and the interior
  of the dual $\op{Conv}(\Lambda_\Gamma)^\ast$. Actually, it was already observed in Thurston's book on the geometry   of topology of three-manifolds (\cite{thurstonbook}) that the action of $\Gamma$ is proper not only on $\HH^n$,   but also in the bigger region $\op{Conv}(\Lambda_\Gamma)^\ast$.
\end{remark}

\subsection{The flat case}
In the flat case, the subgroup $\Gamma$ is not always achronal --- for example, consider the case
of a cocompact lattice in the group of translations of $\RR^{1,n}$. Nevertheless, as in the de Sitter case,
one can show that achronal subgroups are subgroups of Isom$(\RR^{1,n})$ that preserve a closed
achronal subset $\Lambda_\Gamma$ in the Penrose conformal boundary $\partial$Mink$^+(x_0) = C(x_0)$ of a point $x_0$ in $\Ein^{1,n}$. More precisely, such a group always preserves the ``spatial infinity''
$i_0 = x_0$, hence we have to be more precise and distinguish several cases:

-- \textbf{Case $(1)$}: $\Gamma$ is relatively compact: then $\Gamma$ is achronal and no $\Gamma$-orbit accumulates
at the conformal boundary; the limit set $\Lambda_\Gamma$ is then the empty set, and if $\Gamma$ is discrete, the action of $\Gamma$
on $\RR^{1,n} = \Omega(\emptyset)$ is proper.

-- \textbf{Case $(2)$}: $\Gamma$ is not relatively compact, but its linear part is relatively compact in $\SO_0(1,n)$: then, $\Gamma$ preserves
a flat euclidian metric in $\RR^{1,n}$; hence a foliation by parallel spacelike hyperplanes. Every $\Gamma$-orbit is achronal, and accumulates at the spatial infinity $i_0$, that we define to be the limit set. If $\Gamma$ is discrete, the action on the entire Minkowski space, which is the invisible domain for $x_0$, is proper.
Observe that according to the Bieberbach Theorem, if $\Gamma$ is discrete, then, up to finite index, $\Gamma$ is
a group of spacelike translations.

-- \textbf{Case $(3)$}: the linear part of $\Gamma$ is not relatively compact in $\SO_0(1,n)$ and admits a limit set
$\Lambda^0_\Gamma$ in $\SS^{n-1}$. Then, $\Gamma$ is achronal if and only if it preserves an achronal set
in $\mathcal I^+$. The set of accumulation points in $\mathcal I^\pm \approx \SS^{n-1} \times \RR$ is a limit set $\Lambda^\pm_\Gamma$, which is the graph of a map $f:\Lambda^0_\Gamma \to \RR$ (not necessarily Lipschitz). The interior of the set of achronal $\Gamma$-orbits in the set $\Omega(\Lambda_\Gamma)$, which is the union of the two globally hyperbolic domains $\Omega(\Lambda^\pm_\Gamma)$, on which, if it is discrete, $\Gamma$ acts properly.


\subsection{The anti-de Sitter case}
In the anti-de Sitter case, we still have the situation that subgroups of $\SO_0(2,n)$ may fail to be achronal.
The criterion is simpler than in the flat case: $\Gamma$ is achronal if and only if it preserves a closed achronal subset in $\partial\AdS^{1,n}$. Then, one can define a limit set $\Lambda_\Gamma$ as the set of accumulation points
of achronal orbits. This is a closed achronal subset, contained in any $\Gamma$-invariant closed achronal subset of $\partial\AdS^{1,n}$.
The interior of the set of achronal $\Gamma$-orbits in $\AdS^{1,n}$ is the invisible domain $\Omega(\Lambda) = E(\Lambda) \cap \AdS^{1,n}$, on which $\Gamma$ acts properly, even if the quotient is not always globally hyperbolic (it is globally hyperbolic if and only if $\Lambda_\Gamma$ is edgeless, \ie is the graph of a $1$-Lipschitz map defined on the entire sphere $\SS^{n-1}$).

There is a very interesting criterion for achronality, involving a particular element of the bounded cohomology of $\SO_0(2,n)$. Recall the exact sequence:
$$0 \to \ZZ \to \wt\SO_0(2,n) \to \SO_0(2,n) \to 0$$
where the $\ZZ$ is the cyclic group generated by the Galois automorphism $\delta$.
Then, there is a canonical non-algebraic section $s: \SO_0(2,n) \to \wt\SO_0(2,n)$: loosely speaking,
$s(g)$ is the only lift of $g \in \SO_0(2,n)$ such that for every affine domain $U$ in $\wt\AdS^{1,n}$,
the intersection $U \cap s(g)U$ is never empty. Then, for every $g_1$, $g_2$ in $\SO_0(2,n)$, let
$c(g_1, g_2)$ be the integer in $\{ -1, 0, +1\}$ characterized by $s(g_1g_2) = \delta^{c(g_1, g_2)}s(g_1)s(g_2)$.
Then, $c$ is a cocycle, representing an bounded cohomology class in $H^2_b(\SO_0(2,n), \ZZ)$, called the
\emph{bounded Euler class.} \index{bounded Euler class} Then, for every subgroup $\Gamma$ of $\SO_0(2,n)$ the restriction of $c$ to
$\Gamma$ is an element of $H^2_b(\Gamma, \ZZ)$, called the bounded Euler class of $\Gamma$ and denoted by
eu$_b(\Gamma)$.

\begin{teo}[\cite{quasifuchs}]\label{teo:euler}
  A subgroup $\Gamma$ of $\SO_0(2,n)$ preserves a closed achronal subset of $\partial\AdS^{1,n}$ if and only if
  $\op{eu}_b(\Gamma)=0$.
\end{teo}

\subsection{Anosov representations}
\label{sub:anosov}
Recall that for any Lie group $G$, we denote by Rep$(\Gamma, G)$ the modular space of representations of $\Gamma$ in $G$ up to conjugacy in the target $G$.
$$\op{Rep}(\Gamma, G) := \op{Hom}(\Gamma, G)/G.$$
In this section, we restrict to the case where $\Gamma$ is a Gromov hyperbolic group. We denote by $\partial\Gamma$ the Gromov boundary of $\Gamma$, and by $(\wt{U\Gamma,} \tilde{\phi}^t)$ the \emph{geodesic flow of $\Gamma$} \index{geodesic flow} (see \cite{champetier, mineyev}): there is a proper cocompact action of $\Gamma$ on $\wt{U\Gamma}$ commuting with $\tilde{\phi}^t$ such that
the induced flow $\phi^t$ on the quotient space $U\Gamma := \Gamma\backslash\wt{U\Gamma}$ has a hyperbolic behavior. In particular, there are two $\Gamma$-equivariant maps $\xi^\pm:  \wt{U\Gamma} \to \partial\Gamma$ which are constant along the orbits of $\tilde{\phi}^t$; more precisely, one can see $\xi^+(p)$ (respectively $\xi^-(p)$)
as the limit of $\tilde{\phi}^t(p)$ for $t \rightarrow +\infty$ (respectively for $t \rightarrow -\infty$).
There is also a map $\nu^+$ defined on $\wt{U\Gamma}$ that associates with every $p$ a metric $\nu^+(p)$ in a neighborhood of $\xi^+(p)$, such that for every $\gamma$ in $\Gamma$, the action of $\gamma$ on $\partial\Gamma$ is an isometry between the metric $\nu^+(p)$ near $\xi^+(\gamma)$ and the metric $\nu^+(\gamma.p)$ near $\xi^+(\gamma.p)$. Moreover, $\nu^+(\tilde{\phi}^t(p))$ increases exponentially with $t$.

Let $X$ be a manifold on which $G$ acts analytically, typically, a homogeneous manifold $G/H$.
A representation $\rho: \Gamma \to G$ is\emph{ $(G,X)$-Anosov} \index{$(G,X)$-Anosov representation} if the properties of the geodesic flow still hold when we replace $\partial\Gamma$ by $X$. More precisely:

-- there is a continuous equivariant map $f: \partial\Gamma \to X$, where $\partial\Gamma$ is the Gromov boundary of $\Gamma$,

-- there is a continuous $\Gamma$-equivariant family of metrics $\nu^+_\rho(p)$ parameterized by $\wt{U\Gamma}$, so that $\nu^+_\rho(p)$ is a metric on $X$ near $f(\xi^+(p))$ that increases exponentially along the orbits of $\tilde{\phi}^t$.

For more details or other presentations of the notion of Anosov representations, see \cite{guichard3,guichard4, guichard5} or the section $2.1$ of \cite{pressuremetric} or also
the article \cite{canarysurvey}.

The notion of Anosov representations was introduced by Labourie (\cite{labourie}). It has several interesting features:

-- Anosov representations form an open domain $\op{Rep}_{an}(\Gamma, G)$ of $\op{Rep}(\Gamma, G)$,

-- Anosov representations are always faithful with discrete image.

Labourie proved in particular that for any $n \geq 2$, and in the case where $\Gamma$ is the fundamental group of a closed surface, one connected component (the so-called \emph{Hitchin component}) \index{Hitchin component} of $\op{Rep}(\Gamma, \op{PSL}(n,\RR))$ is made of $(\op{PSL}(n, \RR), \mathcal F_n)$-Anosov representations, where $\mathcal F_n$
is the variety of complete flags in $\RR^n$.

\medskip

\noindent\textbf{The de Sitter case.}
Once more, de Sitter space is the dual of hyperbolic space $\HH^{n+1}$. It is well-known (see \cite{guichard3, canarysurvey}) that for any subgroup $\Gamma$  of $\SO_0(1,n+1)$, the inclusion $\Gamma \subset \SO_0(1,n+1)$
is $(\SO_0(1,n+1), \SS^n)$-Anosov if and only if it is convex cocompact.

A particular interesting family of convex cocompact subgroups are the \emph{quasi-Fuchsian subgroups,} \index{quasi-Fuchsian subgroup} \ie discrete subgroups isomorphic to a cocompact lattice of $\SO_0(1,n+1)$ such that the limit set $\Lambda_\Gamma$ is a topologically embedded sphere in $\partial\HH^{n+1}$.

We point out here an interesting characterization of convex cocompact subgroups appearing in our framework:
\emph{a non-elementary discrete subgroup $\Gamma$  of  $\SO_0(1,n+1)$ is convex cocompact if and only the locally de Sitter globally hyperbolic domain
$\Gamma\backslash\Omega_\pm(\Lambda_\Gamma)$ is spatially compact.}

In other words,\emph{ one can characterize $(\SO_0(1,n+1), \SS^n)$-Anosov representations as the holonomy representations of MGHC de Sitter spacetimes that are quotients of hyperbolic domains of $\dS^{1,n}$} (observe that this last condition is restrictive: there are closed M\"obius manifolds whose developing map is not injective; therefore their associated MGHC de Sitter spacetimes have also a non-injective developing map).

\begin{rema}
  Unlike the other cases (see below), $\op{Rep}_{an}(\Gamma, \SO_0(1,n+1))$ is \textbf{not} a closed subset of $\op{Rep}(\Gamma, \SO_0(1,n+1))$. Indeed, quasi-Fuchsian subgroups can be continuously deformed to the trivial subgroup. It follows that the holonomy representation of a MGHC de Sitter spacetime, even if it is Gromov hyperbolic, is not necessarily $(\SO_0(1,n+1), \SS^n)$-Anosov.
\end{rema}

\medskip

\noindent\textbf{The anti-de Sitter case.}
The de Sitter case has been fully treated in \cite{merigot}:

\begin{teo}[Theorem $1.2$ in \cite{merigot}]
  Let $\Gamma$ be the fundamental group of a closed manifold of dimension $n$. Assume that $\Gamma$ is Gromov hyperbolic. Then, a representation
  $\rho: \Gamma \to \SO_0(2,n)$ is $(\SO_0(2,n), \Ein^{1,n-1})$-Anosov if and only if it preserves
  a closed edgeless \textbf{acausal} subset $\Lambda$ of $\Ein^{1,n-1}$.
\end{teo}


\begin{rema}
   If $\Gamma$ is a discrete subgroup acting on $\AdS^{1,n}$ and preserving a closed edgeless acausal subset $\Lambda$ of $\Ein^{1,n-1}$
   such that the quotient $\Gamma\backslash\Omega(\Lambda)$ is spatially compact, then it is Gromov hyperbolic (see Section $8.3.2$ in \cite{merigot}). Therefore, $(\SO_0(2,n), \Ein^{1,n-1})$-Anosov representations are precisely holonomy representations of MGHC $\AdS$-spacetimes whose limit set is acausal, and not simply achronal.
 \end{rema}

 \begin{rema}
 Split $\AdS$-spacetimes have a fundamental group isomorphic to a lattice of $\SO_0(1,k) \times \SO_0(1, \ell)$, hence are not Gromov hyperbolic. We already mentioned that they may be exactly the MGHC $\AdS$-spacetimes with non-acausal limit set.
 \end{rema}

 Finally, we have the following:

 \begin{teo}[\cite{quasifuchs}]
 Let $\Gamma$ be a Gromov hyperbolic group, isomorphic to the fundamental group of a closed $n$-dimensional manifold. Then, if not empty, the space $\op{Rep}_{an}(\Gamma, \SO_0(2,n))$ of Anosov representations is open and closed in the space $\op{Rep}(\Gamma, \SO_0(2,n))$. In particular,  it is a union of components of $\op{Rep}(\Gamma, \SO_0(2,n))$.
 \end{teo}

One of the main intermediate steps in the proof of this theorem is the following result: if $\Gamma \subset \SO_0(2,n)$ is a Gromov hyperbolic group and the holonomy group of a MGHC $\AdS^{1,n}$-spacetime, then any achronal $\Gamma$-invariant closed edgeless achronal subset of $\partial\AdS^{1,n}$ is automatically acausal.

\medskip
\noindent\textbf{The flat case.}
This is the only case that has not been yet completely studied. The only published result is
the work of S. Ghosh (\cite{ghosh}): let $\Gamma$ be a discrete subgroup of Isom$(\RR^{1,2})$, acting
properly discontinuously on $\RR^{1,2}$, and admitting as linear part a convex cocompact subgroup of
$\SO_0(1,2)$. Then, the inclusion $\Gamma \subset$ Isom$(\RR^{1,2})$ is $(\op{Isom}(\RR^{1,2}), \mathcal I^+)$-Anosov. This result is the starting point for the construction of a certain metric on the space of such representations, involving the thermodynamical formalism, the \textit{pressure metric.}

 We propose here the following conjecture: \emph{let $\Gamma$ be a Gromov hyperbolic group. Then a representation
 $\rho: \Gamma \to \op{Isom}(\RR^{1,n})$ is $(\op{Isom}(\RR^{1,n}), \mathcal I^+)$-Anosov if and only if it
 is achronal and the associated globally hyperbolic domain $\Omega(\Lambda)$ is proper.}\footnote{Recall that
 $\Omega(\Lambda)$ is proper if its closure contains no affine line.}

 These representations should also be characterized by the following property: \emph{their linear part $L_\rho: \Gamma \to \SO_0(1,n)$ is faithful and convex cocompact.}

\section{The three-dimensional case: links with Teich\-m\"uller space}
\label{sec:trois}
In this section, we describe the several connections between the Teichm\"uller space and MGHC $\AdS^{1,2}$-spacetimes revealed in Mess's work (\cite{mess1, mess2}) and followers. There are also similar links with MGHC spacetimes of constant curvature $0$ or $+1$,
but we have decided to focus on the $\AdS$ case. For more information on this topic, see \cite{carlip, benbon, SSCH}.

\subsection{Anti-de Sitter space as a space of matrices}
\label{sub:adsmatrix}
\noindent\textbf{An alternative description of $\AdS^{1,2}$.}
Consider the linear space Mat$(2, \RR)$ of two-by-two matrices with real coefficients, equipped with the quadratic form $-$det.
This space is isometric to $(\RR^{2,2}, q_{2,2})$, hence $\AdS^{1,2}$ is naturally identified with the space of matrices of determinant $1$, \ie SL$(2, \RR)$.

We have also the formula $-\det(A) {\rm Id}= A\widehat{A} = \widehat{A}A$ where ${\rm Id}$ is the identity matrix and $\widehat{A}$ the transpose of the matrix of cofactors. Therefore the bilinear form associated with $-\det$ is:
$$\langle X \mid Y \rangle = \op{Tr}(\widehat{X}Y)/2 = \op{Tr}(X\widehat{Y})/2$$
where Tr denotes the trace.
\medskip

\noindent\textbf{The isometry group.}
The actions of SL$(2, \RR)$ on itself by left and right translations are isometric, and they provide a natural identification between $\SO_0(2,2)$ and the quotient of SL$(2, \RR) \times$ SL$(2, \RR)$ by the group of order two generated by $(-{\rm Id}, -{\rm Id})$. Therefore, a representation $\rho: \Gamma \to \SO_0(2,2)$ splits in two representations $\rho_L$, $\rho_R$ from $\Gamma$ into PSL$(2, \RR)$.
\medskip

\noindent\textbf{The Minkowski space and the hyperbolic plane.}
Let $M$ be the $3$-dimensional linear subspace of Mat$(2, \RR)$ consisting of matrices with null trace.
This is the tangent space of SL$(2, \RR)$ at the identity map.
The restriction of $-$det to $M$ has index $(1,2)$, hence $(M, -\op{det})$ is a model for $\RR^{1,2}$. Observe also that for $X$ in $M$, $-$det$(X)$ coincides with Tr$(X^2)/2$ in $M$.
In particular, the hyperbolic plane $\HH^2$ is naturally identified with the space of matrices of determinant $1$ and trace $0$, therefore, with the space of complex structures $J$ on $\RR^2$ (since $J^2 + {\rm Id} = 0$).
In this model, the isometric action of an element $A$ of SL$(2, \RR)$ is the action by conjugacy $J \mapsto AJA^{-1}$.

\medskip

\noindent\textbf{The Einstein spacetime $\Ein^{1,1}$.}
The Klein model of the boundary of $\AdS^{1,2}$ is the projectivization of the $-$det null cone in
 Mat$(2, \RR)$, \ie of non zero non-invertible matrices. Such a matrix determines two lines in $\RR^2$: its image and its kernel; moreover the left action PSL$(2, \RR)$ preserves the kernel and the right action PSL$(2, \RR)$ preserves the image. Therefore, there is a canonical identification of $\overline{\Ein}^{1,1}$ with $\RR\mathbb P^1 \times \RR\mathbb P^1$, where
 the left factor of PSL$(2, \RR) \times$ PSL$(2, \RR)$ acts trivially on the left factor $\RR\mathbb P^1$, and the right factor acts trivially on the right factor $\RR\mathbb P^1$. The lightlike geodesics of $\overline{\Ein}^{1,2}$ are the factors $\{ \ast \} \times \RR\mathbb P^1$ and the factors $\RR\mathbb P^1 \times \{ \ast \}$. The image of the isotropic cone $\mathbb{P}(\mathcal{C}_2))$ in $\mathbb{P}(\RR^{2,2})$ is also ruled by two families of lines: they are these two families of lightlike geodesics.

 This description of $\overline{\Ein}^{1,1}$ as a product of two circles is different from the analogous decomposition for its double cover $\Ein^{1,1}$, which was defined as the product
 $\SS^1 \times \SS^1$ equipped with the metric $\bar{g}_1 - \bar{g}_1$. The factors $\{ \ast \} \times \RR\mathbb P^1$ and $\RR\mathbb P^1 \times \{ \ast \}$ can be seen, locally, as the diagonal and anti-diagonal in $\SS^1 \times \SS^1$.
 \medskip

 \noindent\textbf{Achronal subsets of $\Ein^{1,1}$.}
Every closed acausal subset of $\overline{\Ein}^{1,1}$ intersects every circle $\RR\mathbb P^1 \times \{ \ast \}$ and $\{ \ast \} \times \RR\mathbb P^1$ in at most one point. Therefore, these sets are graphs of increasing maps from $\RR\mathbb P^1$ into $\RR\mathbb P^1$. In particular, closed edgeless acausal subsets are precisely the graphs of the homeomorphisms of $\RR\mathbb P^1$.

At the limit, every achronal subset $\Lambda$ is a \emph{generalized graph of a semi-conjugacy.} By this we mean
that $\Lambda$ is the union of the graph of a non-decreasing map $f: \Lambda_0 \to \RR\mathbb P^1$ ($\Lambda_0$ is a closed subset of $\RR\mathbb P^1$) and some vertical segments $\{ \ast \} \times I$ where $I$ is a closed segment in $\RR\mathbb P^1$ (there is also the limiting case of purely lightlike subsets, where $\Lambda_0 =\emptyset$ and where
$\Lambda$ is one factor $\{ \ast \} \times \RR\mathbb P^1$).

We hope that the reader will easily agree with the idea that the less awkward definition of ``generalized graph of generalized maps from $\RR\mathbb P^1$ onto $\RR\mathbb P^1$''
is the first definition given here, as closed achronal subsets of $\overline{\Ein}^{1,1}$.
 \medskip

 \noindent\textbf{Timelike geodesics.}
 We invite the reader to keep in mind the content of the end of Section \ref{sub:adsbasic}. See also the reference \cite[Section $7.3$]{barbtz1}.
In this low dimension, the space of timelike geodesics $\cT_{4} := \SO_0(2,2)/(\SO(2) \times \SO(2)) \approx (\op{PSL}(2,\RR) \times \op{PSL}(2,\RR))/(\SO(2) \times \SO(2))$
is isometric (up to a constant factor) to $\HH^2 \times \HH^2$ equipped with the product metric.
This correspondence is expressed in the following way: the timelike geodesic corresponding to an element $(x,y)$ of
$\HH^2 \times \HH^2$ is the unique timelike geodesic preserved by the subgroup of PSL$(2, \RR) \times$ PSL$(2, \RR)$ fixing $(x,y)$, namely made of pairs $(g, h)$ where $g$ is a rotation at $x$ and $h$ a rotation at $y$.
In particular, for every $g$ in PSL$(2, \RR) \approx \AdS^{1,2}$, the set of timelike geodesics containing $g$ is the graph of $g$ in $\HH^2 \times \HH^2$.

The symplectic form $\omega$ is $\omega_0 - \omega_0$ where $\omega_0$ is the volume form on $\HH^2$.

The following observation will be useful later:
Let $(x_1, y_1)$ and $(x_2,y_2)$ be two elements of $\HH^2 \times \HH^2$.
Then, the associated timelike geodesics of $\AdS^{1,2}$ have a common point if and only if the hyperbolic distance between $x_1$ and $x_2$ is equal to the distance between $y_1$ and $y_2$ --- common points are then the elements of
PSL$(2, \RR)$ mapping $x_1$ to $y_1$ and $x_2$ to $y_2$.
\medskip

\noindent\textbf{Timelike geodesics 2: more explicit computations.}
Recall that $\mathcal U^{1,2}$ is the space of future oriented timelike vectors of $\AdS^{1,2}$.
In our matrix model, $\mathcal U^{1,2}$ is the space of pairs of matrices $(A, C)$ satisfying:

-- det$(A) =$ det$(C) = 1$,

-- $A$ and $C$ are orthogonal for $-$det.

The last condition is equivalent to Tr$(AC^{-1}) = 0$. But recall that the space of matrices of trace zero and determinant $1$ is a model
for $\HH^2$.
Therefore, $(A,C) \mapsto (J_L, J_R)$, where $J_L := AC^{-1}$ and $J_R := C^{-1}A$, defines a map between $\mathcal U^{1,2}$ and $\HH^2 \times \HH^2$. This map is constant along the geodesic flow, and induces (up to a scalar constant) the $\SO_0(2,2)$-equivariant isometry between $\cT_4$ and $\HH^2 \times \HH^2$ - it is clear that $J_L$ (respectively $J_R$) is invariant by right translations (respectively left translations).

\subsection{\mathversion{bold}$\AdS^{1,2}$ globally hyperbolic spacetimes}\label{sub:ads3}
References adapted to the content of this section are \cite{benbon}  and more recently \cite{SSCH}.
\medskip

\noindent\textbf{Mess's parametrization by Teich\mathversion{bold}$(\Sigma) \times$Teich\mathversion{bold}$(\Sigma)$.}
Let $M_\Lambda(\Gamma) = \Gamma\backslash\Omega(\Lambda)$ be a MGHC spacetime locally modeled on $\AdS^{1,2}$.
Then $\Lambda$ is the generalized graph of a semi-conjugacy between the projection $\Gamma_L$ of $\Gamma$ in the left factor PSL$(2,\RR)$ and the projection $\Gamma_R$ of the projection in the right factor. It follows that
these projections are injective and that $\Gamma_L$ and $\Gamma_R$ are both discrete subgroups of PSL$(2,\RR)$.
Moreover, these groups are isomorphic to the fundamental group of Cauchy surfaces of $M_\Lambda(\Gamma)$, hence
of a closed surface. It follows that $\Gamma_L$ and $\Gamma_R$ are cocompact lattices, and that
$\Lambda$ is the graph of the unique homeomorphism conjugating $\Gamma_L$ and $\Gamma_R$ --- in particular, it is acausal.

Conversely, given any pair $\rho_L$, $\rho_R$ of faithful and discrete representations into PSL$(2,\RR)$ of the fundamental group $\Gamma$ of a closed surface, the quotient of $\Omega(\Lambda)$ (where $\Lambda$ is the graph of the unique conjugacy between $\rho_L$ and $\rho_R$) by the image of the representation $(\rho_L, \rho_R)$ from
$\Gamma$ into $(\op{PSL}(2,\RR) \times \op{PSL}(2,\RR))/(-{\rm Id} , -{\rm Id}) \approx \SO_0(2,2)$ is MGHC.

Therefore, for every closed surface $\Sigma$, there is a canonical one-to-one correspondence between the space of $\AdS$ MGHC spacetimes diffeomorphic to $\Sigma \times \RR$ up to isometry, and the product Teich$(\Sigma) \times \op{Teich}(\Sigma)$.
\medskip

\noindent\textbf{Diallo's parametrization by Teich\mathversion{bold}$(\Sigma) \times$Teich\mathversion{bold}$(\Sigma)$.}
The boundary of the convex hull $\op{Conv}(\Lambda)$ in $\AdS^{1,2}$ is the union of two spacelike surfaces $S^\pm(\Lambda)$ (its past and future components). Even if they are not smooth, the metric induced on each of them is isometric to $\HH^2$, therefore the quotient surfaces $\Gamma\backslash S^\pm(\Lambda)$ represents two hyperbolic closed surfaces, \ie a point in Teich$(\Sigma) \times \op{Teich}(\Sigma)$ (not to be confused with
the Mess parameters). In his Ph D. thesis (\cite{boubacar}) B. Diallo proved that any element of Teich$(\Sigma) \times \op{Teich}(\Sigma)$ can be obtained in this way, but the uniqueness of the $\AdS^{1,2}$ spacetime realizing this
pair of metrics is still an open question. This result is the Lorentzian analogue of Epstein-Marden's Theorem establishing the realization of any element of Teich$(\Sigma) \times \op{Teich}(\Sigma)$ as the metric induced
on the boundary of the convex core of quasi-Fuchsian hyperbolic $3$-manifolds (\cite{epsteinmarden}).
\medskip

\noindent\textbf{Mess parametrization by measured geodesic laminations.}
As in the classical situation of quasi-Fuchsian hyperbolic manifolds, the future and past components $S^\pm(\Lambda)$ are \emph{pleated} surfaces \index{pleated surface}, embedded isometric copies of $\HH^2$ bended
in $\AdS^{1,2}$ along some
measured geodesic lamination $\lambda^\pm$. The pair of laminations $(\lambda^+, \lambda^-)$ is \emph{filling} \ie
any closed curve $c$ in S which is not homotopically trivial has non-zero intersection with either $\lambda^+$ or $\lambda^-$. Mess proved that the map associating with every $\AdS^{1,2}$ MGHC spacetime the
measured geodesic lamination $(\Gamma\backslash S^+(\Lambda), \lambda^+)$ (or  $(\Gamma\backslash S^-(\Lambda), \lambda^-)$) realizes a one-to-one correspondence.

One can also forget the hyperbolic metric on $\Gamma\backslash S^\pm(\Lambda)$. In \cite{bonschlenk}
Bonsante and Schlenker proved that any filling pair of measured laminations $(\lambda^+, \lambda^-)$ (but forgetting the hyperbolic metrics) is realized as a pair of pleated laminations on the convex core of a MGHC $\AdS$-spacetime.
\medskip

\noindent\textbf{\mathversion{bold}$\AdS^{1,2}$-spacetimes and earthquakes.}
In \cite{mess1} Mess pointed out a connection between $\AdS$ spacetimes and the notion of \emph{earthquake} \index{earthquake} introduced by Thurston (see \cite{kerk}). Let $\lambda$ be a measured geodesic lamination on $\HH^2$. A convenient way to see $\lambda$ is to see it as a stratification of $\HH^2$, with one-dimensional strata (the geodesics in $\lambda$) and $2$-dimensional strata (the closure of the components of $\HH^2 \setminus \op{Supp}(\lambda)$, where Supp$(\lambda)$ is the support of $\lambda$). The (left) earthquake defined by $\lambda$ is a non-continuous map $E_\lambda: \HH^2 \to \HH^2$ that is an isometry on each $2$-dimensional stratum, but for which each $1$-dimensional stratum is a \emph{rift} on which $E_\lambda$ may be non continuous. When $\lambda$ is a rational lamination, \ie a locally finite collection of weighted geodesic, we can give a more precise definition: if $C_1$ and $C_2$ are two $2$-dimensional strata bounding the same leaf $\ell$ of $\lambda$, and if $g_1$ and $g_2$ are the isometries of $\HH^2$ coinciding with $E_\lambda$ on $C_1$ and $C_2$, respectively, then $g_2 = h_a \circ g_1$ where $h_a$ is the element of PSL$(2,\RR)$ defined as follows:

-- orient the geodesic $\ell$ so that $C_1$ is on the left of $\ell$, and $C_2$ on the right of $\ell$,

-- let $a$ be the weight of $\ell$ for $\lambda$,

then $h_a$ is the unique hyperbolic element preserving $\ell$ so that for every $x$ in $\ell$ the image $h_a(x)$ is the unique element of $\ell$ in the direction of the orientation of $\ell$ and at distance $a$. The right earthquake is obtained with the other orientation of $\ell$.

Earthquakes for general measured laminations are then defined by a limiting process, at least in the case of
measured geodesic laminations invariant by a cocompact lattice of PSL$(2, \RR)$, involving the density of rational laminations.

Earthquakes have then a natural extension to the boundary $\partial\HH^2$. When $\lambda$ is preserved by
a cocompact lattice, this extension is a homeomorphism, but if not, this extension in general is merely a ``generalized semi-conjugacy'', \ie in the point of view adopted in this survey, a closed edgeless achronal subset of $\partial\HH^2 \times \partial\HH^2$.
A celebrated theorem by Thurston is that\emph{ any homeomorphism of the circle is realized by a left earthquake,}
and our purpose here is to point out that this theorem is well explained in our $\AdS$-background. Moreover, our presentation has the advantage of providing a simpler and more direct definition of earthquakes avoiding the technical difficulties associated with the step between rational measured laminations and general ones.

Indeed, let $\Lambda$ be a closed achronal edgeless subset of $\Ein^{1,1}$ (for example, the graph of a homeomorphism). Consider the globally hyperbolic domain $E(\Lambda)$. Then, the space $\Sigma(\Lambda)$ of cosmological lines is a closed embedded locally Lipschitz disk in $\cT_4 \approx \HH^2 \times \HH^2$.
This is not the graph of a map from $\HH^2$ into $\HH^2$ but almost: this is directly the left earthquake. Indeed,
this is the set of timelike geodesics orthogonal to support hyperplanes of the future component $S^+(\Lambda)$ of the convex hull of $\Lambda$, \ie the level set $\{ \tau = \pi/2 \}$ of the cosmological time. The geodesics orthogonal to a totally geodesic face of the pleated surface $S^+(\Lambda)$ form a region of $\Sigma(\Lambda)$ which is the graph of
the restriction of an element of PSL$(2, \RR)$ to a region of $\HH^2$ bounded by geodesics of $\HH^2$, whereas the timelike geodesics orthogonal to $S^+(\Lambda)$ at the bending locus $\lambda^+$ form some ``vertical bands'' made of points $(x, y)$ where $x$ describes a geodesic $\ell_1$ of $\HH^2$ and $y$ some segment $I$ in a geodesic $\ell_2$ of $\HH^2$. The hyperbolic length of $I$ is the measure of the associated leaf of $\lambda^+$.

In summary, in this vision, a left earthquake is the set of timelike geodesics orthogonal to a pleated surface, and the main idea of Thurston's Earthquake Theorem reduces to the fact that the boundary of the convex hull of a closed edgeless achronal subset is a pleated surface.
\medskip

\noindent\textbf{Cauchy surfaces and volume preserving maps.}
Anti-de Sitter geometry also proposes a new vision on area preserving maps between hyperbolic surfaces.
Indeed, let $M_\Lambda(\Gamma)$ be a MGHC $\AdS^{1,2}$ spacetime, defined by two Fuchsian representations $\rho_L, \rho_R: \Gamma \to$ PSL$(2,\RR)$. For every Cauchy surface $S$ of $M_\Lambda(\Gamma)$, the Gauss map of $S$ provides a $\Gamma$-equivariant map $\nu: \tilde{S} \to \HH^2 \times \HH^2$ whose image is $\omega$-Lagrangian.
Since $\omega = \omega_0 - \omega_0$, in the region where $\nu(\tilde{S})$ is locally the graph of a map $f$ from $\HH^2$ into $\HH^2$, $f$ is a volume preserving map.

We can be actually more precise: recall that the map $\mathcal U^{1,2} \to \HH^2 \times \HH^2$ is given, in term
of two-by-two matrices, by $(A, C) \mapsto (J_L = AC^{-1}, J_R = C^{-1}A)$.
The derivative of this map at a point $(A,C)$ is given by:
\begin{eqnarray*}
  \dot{J}_L &=& \dot{A}C^{-1} - AC^{-1}\dot{C}C^{-1} = [\dot{A} - J_L\dot{C}]C^{-1} \\
  \dot{J}_R &=& C^{-1}[\dot{A} - \dot{C}J_R].
\end{eqnarray*}
Hence the hyperbolic norm $-\op{det}(\dot{J}_L)$ is $-\det(\dot{A} - J_L\dot{C})$ and
the norm of $-\op{det}(\dot{J}_L)$ is $-\det(\dot{A} - \dot{C}J_R)$.

Now let $p \mapsto (A(p), C(p))$ be the spacelike immersion induced by a smooth spacelike immersion from a surface $\tilde{S}$ into $\AdS^{1,2}$: $A(p)^\perp \cap C(p)^\perp$ is the plane tangent to $A(S)$ at $A(p)$. Then:

-- The map $X \mapsto J_LX$ preserves $A(p)^\perp \cap C(p)^\perp$: for $X \in A(p)^\perp \cap C(p)^\perp$, we have:
$$\op{Tr}(\widehat{A}J_LX) = \op{Tr}(A^{-1}AC^{-1}X)= \op{Tr}(C^{-1}X)=0,$$
and, since $J_L^2 = -I_2$:
$$\op{Tr}(\widehat{C}J_LX) = \op{Tr}(C^{-1}AC^{-1}X)= \op{Tr}(A^{-1}AC^{-1}AC^{-1}X = \op{Tr}(A^{-1}X)=0.$$

-- For $X \in A(p)^\perp \cap C(p)^\perp$ we have:
$$\langle J_LX \mid X \rangle = \op{Tr}(J_LX\widehat{X})= \det(X)\op{Tr}(J_L)=0.$$
It follows that $J_L$ is the complex structure on the tangent space associated
with the metric $I$ on $\tilde{S}$ induced by the $\AdS$-metric.

-- Similar computations show that $X \mapsto -XJ_R$ is also this complex structure, therefore:
$$\forall X \in A(p)^\perp \cap C(p)^\perp \;\;XJ_R = -J_LX.$$

Therefore:
\begin{eqnarray*}
  \dot{J}_L &=&  [\dot{A} - J_L\dot{C}]C^{-1} \\
  \dot{J}_R &=& C^{-1}[\dot{A} + J_L \dot{C}].
\end{eqnarray*}

We deduce the formula relating the metric $I$, the associated complex structure $J_L$, the shape operator $B$ and the left or right hyperbolic metrics $g_L$, $g_R$ induced by the composition of the Gauss map with the projection $\cT_4 \approx \HH^2 \times \HH^2$ on the left factor (compare with \cite[Lemma $2.9$]{SSCH}):
\begin{eqnarray*}
  g_L(v) &=& I(v - J_LB(v))  \\
  g_R(v) &=& I(v + J_LB(v)).
\end{eqnarray*}

In particular, $\nu(\tilde{S})$ fails to be transverse to the fiber of the left or right factor if and only if
there is a non zero tangent vector $v$ at $p$ satisfying $B(v) = \pm J_L(v)$. If this happens, then, since $B$ is self-adjoint and this basis is orthogonal, the symmetric matrix expressing $B$ in the basis $(v, J_L(v))$ is:
$$\pm\left(
  \begin{array}{cc}
    0 & 1 \\
    1 & 0 \\
  \end{array}
\right).$$
In other words, $\nu(\tilde{S})$  is not locally the graph of a map from $\HH^2$ into $\HH^2$ near $p \in \tilde{S}$ if and only its mean curvature at $p$ vanishes and its scalar curvature is $-1$. This condition can also be formulated in another way: it means that $\tilde{S}$ is tangent at order two to the space of timelike geodesics orthogonal to a given spacelike geodesic (compare with the previous study of Gauss maps of pleated surfaces).

In particular, if the shape operator $B$ satisfies I$(B(v)) <$ I$(v)$ everywhere, \ie if the absolute values of the principal curvatures
of $\tilde{S}$ are all $<1$, then $\nu(\tilde{S})$ is the graph of a volume preserving function $f: \HH^2 \to \HH^2$.

Conversely, let $f: \rho_L(\Gamma)\backslash\HH^2 \to \rho_R(\Gamma)\backslash\HH^2$ be any volume preserving smooth map, lifting to a conjugacy $\tilde{f}: \HH^2 \to \HH^2$ between $\rho_L$ and $\rho_R$. The graph of $\tilde{f}$ is a $\Gamma$-invariant embedded Lagrangian submanifold. At the end of Section \ref{sec:models} we have seen that $\tilde{f}$ provides a map $\varphi: \HH^2 \to \mathcal U^{1,2}$, where $\mathcal U^{1,2}$ is the space of future oriented timelike vectors of $\AdS^{1,2}$. We claim that furthermore, this map is $\Gamma$-equivariant. This is not completely obvious since what is immediate is that for every $\gamma$ in $\Gamma$, we have $\varphi \circ \gamma = \Phi^{t(\gamma)} \circ \varphi$ for some real number $t(\gamma)$ which \emph{a priori} might be non-trivial. But $t: \Gamma \to \RR$ is then a homomorphism, and if it is not trivial, then it would have arbitrary big values, and this would contradict the fact that $\Gamma$ preserves a closed edgeless achronal subset of $\wt\AdS^{1,n}$.

Therefore,
the composition of $\varphi$ with the bundle map $\mathcal U^{1,2} \to \AdS^{1,2}$ is a $\Gamma$-equivariant map
from $\HH^2$ into $\AdS^{1,2}$, which is ``almost spacelike'' as explained at the end of Section \ref{sec:models}.

For a more complete treatment (but with a different presentation) of the link between Cauchy surfaces in MGHC $\AdS^{1,2}$-spacetimes and volume preserving maps between hyperbolic surfaces, with a list of special types of spacelike surfaces with remarkable associated volume preserving maps between hyperbolic surfaces, see \cite{bonschlenkmax, SSCH} and also \cite{fillgraham} to appear in this {\it Handbook of Group actions.}

\subsection{BTZ-multi black holes}
\label{sub:btz2}
In this section, we summarize the content of \cite{barbtz1, barbtz2}, to which we refer for further details (see also \cite{ammi1, ammi2, ammi3, ban1, ban2, brill1, brill2} in the physics literature). Consider the case of domains $\Omega(\Lambda)$ of $\AdS^{1,2}$ associated with closed achronal subsets $\Lambda$ of $\Ein^{1,1}$, but not necessarily edgeless. Then, the invisible domain $E_0(\Lambda)$ in
$\Ein^{1,1}$ can be interpreted as the conformal boundary of $\Omega(\Lambda)$. When $\Lambda$ is preserved by
a discrete subgroup $\Gamma$ of $\SO_0(2,2)$, the quotient of $E_0(\Lambda)$ by $\Gamma$ is a conformal boundary for $M_\Lambda(\Gamma)$.

In this case, $\Gamma$ is the image of a representation $(\rho_L, \rho_R): \Gamma \to \op{PSL}(2,\RR) \times \op{PSL}(2,\RR) \approx \SO_0(2,2)$, and $\Lambda$ is the (generalized) graph of a (semi\nobreakdash-)conjugacy between $\rho_L$ and $\rho_R$.

Let $\Lambda^+$, $\Lambda^-$ be the future (respectively past) boundary of $E_0(\Lambda)$. Then, $\Omega(\Lambda^+)$ and $\Omega(\Lambda^-)$ are globally hyperbolic domains in $\Omega(\Lambda)$ that can be respectively interpreted as \emph{(multi\nobreakdash)black holes} and \emph{(multi\nobreakdash)white holes.} Indeed, $\Omega(\Lambda^+)$ is the region in
$\Omega(\Lambda)$ which is invisible from the conformal boundary where the observer is assumed to be located.

Despite the fact that it has not been done anywhere in the literature (as far as we know), this vision can be extended in
higher dimensions, where an $n+1$-dimensional AdS multi black hole could be defined as a locally $\AdS^{1,n}$-spacetime admitting a compact convex core, and as conformal boundary a finite union of MGHC conformally flat spacetimes.
Presumably, this theory coincides with the theory of convex cocompact subgroups of $\SO_0(2,n)$.
Such a theory would be a natural analogue of the theory of Schottky groups,
where convex cocompact means $(\SO_0(2,n), \Ein^{1,n-1})$-Anosov.

\section{Proper actions}
\label{sec:propre}
In this last section, we review what is known about discrete isometry groups acting properly on the entire model spacetime.
We will briefly mention results in this direction, once more focusing on the connection with the space of timelike geodesics.

\subsection{Cocompact actions}
This section is essentially a concise extract from a previous survey \cite{surveyBZ}.
One of the most important results is the completeness of Lorentzian manifolds of constant curvature.
The following result is much more difficult to prove than in the well-known Riemannian case:

\begin{teo}[\cite{carriere, klingler}]
 Every closed Lorentzian manifold of constant curvature is geodesically complete.
\end{teo}

This theorem implies that closed Lorentzian manifolds are quotients by discrete groups of isometries of
the simply connected model spacetimes $\RR^{1,n}$, $\dS^{1,n}$ or $\wt\AdS^{1,n}$.

The classification of closed Lorentzian manifolds of constant curvature therefore reduces to the classification
of groups acting properly and cocompactly. It is essentially solved, except in the $\AdS$ case.
\medskip

\noindent\textbf{The de Sitter case.} The de Sitter case is essentially trivial, due to the \emph{Calabi-Markus phenomenom}: \index{Calabi-Markus phenomenom} a group acting properly discontinuously on $\dS^{1,n}$ is necessarily finite. Indeed, let $S$ be the umbilical sphere, that is the intersection between $\dS^{1,n} \subset \RR^{1,n+1}$ and a spacelike hyperplane $H$ of $\RR^{1,n}$. Then for any $g$ in $\SO_0(1,n+1)$, the iterate $gH$ is a hyperplane, hence it intersects non-trivially
$H$, and $H \cap gH \cap \dS^{1,n}$ is non-empty.

Therefore, since a finite group cannot act cocompactly on $\dS^{1,n}$, there is no closed Lorentzian manifold of
positive constant curvature.
\medskip

\noindent\textbf{The flat case.}

\begin{teo}[\cite{FGH, GKam, GrM}]
Let $M = \Gamma\backslash\RR^{1,n}$ be a closed Lorentzian
flat manifold. Then, up to finite covers, $\Gamma$ is a lattice in a solvable subgroup $G$ of Isom$(\RR^{1,n})$ acting
simply transitively on $\RR^{1,n}$.
\end{teo}

There are many possibilities for the solvable Lie group $G$, see
[43] for concrete constructions in the case of the $3$-dimensional Heisenberg
and SOL groups, and \cite{GrM, Guediri} for a general study.
\medskip

\noindent\textbf{The anti-de Sitter case \mathversion{bold}$\AdS^{1,n}$ for \mathversion{bold}$n\geq 3$.}
According to the Chern-Gauss-Bonnet formula (\cite{chern}), for even-dimensional anti de Sitter manifolds, the
Euler number equals the volume, up to a non-trivial multiplicative constant.
But any compact Lorentz manifold has a vanishing Euler number, since it
possesses a direction field. Therefore, $\AdS^{1,n}$ may admit a compact quotient only when $n$ is even.

Conversely, for any even integer $n=2d$, we have a natural inclusion of $U(1,d)$ in $\SO_0(2,2d)$, and
$U(1,d)$ acts transitively and properly discontinuously on $\AdS^{1,2d}$. Therefore, for any cocompact lattice
$\Gamma$ of $U(1,d)$, the quotient $\Gamma\backslash\AdS^{1,2d}$ is a compact Lorentzian manifold of
dimension $1+2d$. In \cite{zeghib} Zeghib conjectured that any closed Lorentzian manifold of constant negative curvature and of odd dimension $\geq 5$
is of this form (observe that by the inclusion $\Gamma \subset U(1,d) \subset \SO_0(2,2d)$
is rigid (\cite{ragh, weil}), \ie any deformation of $\Gamma$ in $\SO_0(2,2d)$
is contained in a conjugate of $U(1,d)$ in $\SO_0(2,2d)$).
\medskip

\noindent\textbf{The  case of \mathversion{bold}$\AdS^{1,2}$.}
Lattices of $U(1,1) \approx \SO_0(2,2)$ still provide compact quotients of $\AdS^{1,2}$, but there are other
examples. This was first observed by Goldman (\cite{gold}). Kulkarni and Raymond proved the following important structure theorem: let $\Gamma$
be a subgroup of PSL$(2, \RR) \times$ PSL$(2, \RR)$ acting properly discontinuously on PSL$(2, \RR)$. Then, up to finite coverings, $\Gamma$ is isomorphic to the fundamental group $\Gamma_g$ of a closed surface, and more precisely, the image of some faithful representation $(\rho, r): \Gamma_g \to \op{PSL}(2, \RR) \times \op{PSL}(2, \RR)$, where (up to swapping the factors) $\rho$ is a Fuchsian representation.
In his thesis Salein (\cite{salein}) studied these representations and found some new examples, based on the following criterion: a representation  $(\rho, r): \Gamma_g \to \op{PSL}(2, \RR) \times \op{PSL}(2, \RR)$
where $\rho$ is Fuchsian acts properly discontinuously on PSL$(2, \RR)$ if and only if there is a \textbf{$1$-contracting map} $f:\HH^2 \to \HH^2$ such that:
$$\forall \gamma \in \Gamma_g \;\; f \circ \rho(\gamma) = r(\gamma) \circ f.$$
In her thesis (\cite{fannythese, fanny}) F. Kassel proved that this criterion is equivalent to the \emph{strict domination} of $r$ by $\rho$, \ie to the requirement that for every $\gamma$ in $\Gamma_g$, the translation length of $r(g)$ is strictly dominated by the translation length of $\rho(\gamma)$ (with the convention that elliptic and parabolic elements have zero translation length).
Recently, Gu\'eritaud, Kassel and Wolf proved that every connected component of
$\op{Rep}(\Gamma_g, \op{PSL}(2, \RR))$, except the two Teichm\"uller components, contains an element $r$ strictly dominated by some Fuchsian representation $\rho$, providing many new examples (\cite{GKW}, see also \cite{derointholozan, tholozan}).

\begin{rema}
  Salein's criterion has a very nice interpretation in terms of timelike geodesics, pointed out in \cite{DGK}, and that we can state in the following way: \emph{there is a $\Gamma$-equivariant $1$-contracting map $f: \HH^2 \to \HH^2$ if and only there is a $\Gamma$-invariant foliation of $\AdS^{1,2}$ by timelike geodesics.} Indeed, the graph of $f$ in $\HH^2 \times \HH^2 \approx \cT_4$ is a family of timelike geodesics such that for every $(x_1, y_1)$, $(x_2, y_2)$ in this family we have $d(y_1, y_2) \neq d(x_1, x_2)$. We have seen at the end of Section \ref{sub:adsmatrix} that this means precisely that the corresponding timelike geodesics are pairwise disjoint.

  Observe that this criterion is very different from the one concerning Gauss maps of Cauchy surfaces: a contracting map $f$ cannot preserve the volume, which was the condition to be the Gauss image of a Cauchy surface.
\end{rema}

\begin{remark}\label{rk:zeghibfoliate}
  The preceding remark may very well have an extension in the higher dimensional case. It seems related to Zeghib's conjecture about compact quotients of $\SO_0(2,2d)$: indeed, $U(1,d)$ can be characterized as the subgroup of $\SO_0(2,2d)$ preserving a foliation of $\AdS^{1,2d}$ by special timelike geodesics: the intersections between $\AdS^{1,2d}$ and $J$-complex lines in $\RR^{2,2d} \approx \mathbb C^{1,d}$ for some complex structure $J$ calibrated with $q_{2,2d}$. Hence it sounds reasonable to split Zeghib's conjectures in two, let us say, half-conjectures:

  -- a discrete subgroup $\Gamma$ of $\SO_0(2,n)$ acts properly discontinuously and cocompactly on $\AdS^{1,2d}$ if and only if it preserves a foliation by timelike geodesics;

  -- if a discrete subgroup $\Gamma$ of $\SO_0(2,n)$ preserves a foliation by timelike geodesics, then it preserves the foliation by timelike geodesics associated with a complex structure calibrated with $q_{2,2d}$.
\end{remark}

\subsection{Margulis spacetimes}
For a very nice recent review of most of this section (in French),
see Schlenker's text for the S\'eminaire Bourbaki (\cite{schbourbaki}).

As mentioned above, compact complete flat manifolds have (virtually) solvable fundamental groups.
In \cite{milnor} Milnor asked the following question: \emph{does the free
group admit a proper action on $\RR^{1,2}$}?  In \cite{margulis1, margulis2}, following a very partial
hint in \cite{milnor}, Margulis gave a positive answer to this question. Since then,
quotients of $\RR^{1,2}$ by torsion-free discrete groups of isometries are called \emph{Margulis spacetimes.}
\index{Margulis spacetime}

Afterwards, Drumm introduced the notion of \emph{crooked planes,} giving
a more intuitive geometric vision on these spacetimes (\cite{drumm1}), and extended
considerably the list of Margulis spacetimes by proving that every discrete
free subgroup of $\SO_0(1, 2)$ is the linear part of the holonomy of a Margulis
spacetime (\cite{drumm2}).

In \cite{GLM}, Goldman, Labourie and Margulis provided a necessary and sufficient
criterion on the so-called Margulis invariant for the action of a discrete subgroup of
Isom$(\RR^{1,2})$ to be properly discontinuous.

In their remarkable recent work \cite{DGK, DGK2} Danciger, Gu\'eritaud and Kassel elucidated the most important remaining questions on Margulis spacetimes, at least in the case where the linear part of the group is
convex cocompact (\ie in the case where the inclusion $\Gamma \subset \op{Isom}(\RR^{1,2})$ is Anosov,
even if this remark has no fundamental importance in their work):

-- these spacetimes are foliated by timelike geodesics, in particular, they are diffeomorphic to
the interior of a handlebody,

-- they all admit a fundamental domain delimited by crooked planes,

-- they can be seen as infinitesimal versions of \emph{Margulis $\AdS$-spacetimes} in a very precise geometric way.

Concerning the last item, let us note that Margulis $\AdS$-spacetimes are quotients of $\AdS^{1,2}$ by discrete subgroups of $\op{PSL}(2, \RR) \times \op{PSL}(2, \RR)$, images of faithful representations $(\rho, r): \mathbb F \to \op{PSL}(2, \RR) \times \op{PSL}(2, \RR)$, where $\mathbb F$ is a free group and $\rho: \mathbb F \to \op{PSL}(2, \RR)$ a convex cocompact representation --- hence, once more, $(\rho, r)$ is $(\SO_0(2,2), \Ein^{1,1})$-Anosov. Let us mention the last paper (\cite{DGK3}) where the authors prove that, unlike flat Margulis spacetimes, Margulis AdS spacetimes do not necessarily admit ``crooked fundamental domains''.

Finally, for extensions of these results to proper affine actions on flat affine spaces in higher dimensions (but then escaping from the Lorentzian framework), see \cite{smilga1, smilga2, smilga3}.



\section*{Acknowledgments}
This work was partially supported by the ANR project ANR 2011 BS01 003 02 GR ANALYSIS GEOMETRY.
A special thank to F. Fillastre for several useful conversations during the preparation of this manuscript, and another
acknowledgement to the referee who contributed to deeply improve the quality of this text.


\end{document}